\documentclass[5p,twocolumn,10pt,times]{elsarticle}
\pdfoutput=1
\usepackage{amssymb,amsthm}
\usepackage{lineno,amsthm,amsmath,mathrsfs,bm}
\usepackage{algorithm,algorithmicx,algpseudocode}
\usepackage{subfigure}
\usepackage{booktabs,multirow,threeparttable}
\usepackage{caption}
\usepackage{float}
\usepackage{graphicx}
\usepackage{epstopdf}
\usepackage{hyperref}
\usepackage{color}

\begin{document}
\baselineskip11pt
\captionsetup[figure]{labelfont={bf},name={Fig.},labelsep=period}

\begin{frontmatter}

\title{Collocation and Mass Matrix in Least-squares Isogeometric Analysis}

\author[zju]{Gengchen Li}
\ead{gengchenli@zju.edu.cn}
\author[zju]{Hongwei Lin\corref{cor}}
\ead{hwlin@zju.edu.cn}
\cortext[cor]{Corresponding author.}
\address[zju]{School of Mathematical Sciences, Zhejiang University, Hangzhou, 310058, China}

\begin{abstract} 
In this paper, we conduct a systematic numerical analysis of the spectral properties of
      the collocation and mass matrices in the isogeometric 
      least-squares collocation method (IGA-L), 
      for the approximation of the Poisson problem with 
      homogeneous Dirichlet boundary conditions. 
This study primarily focuses on the spectral properties of the IGA-L  
	collocation and mass matrices in relation to the isogeometric discretization 
	parameters, such as the mesh size, degree, regularity, spatial dimension, 
      and the number and distribution of the collocation points.
Through a comprehensive numerical investigation, we provide estimations for the  
	condition number, as well as the maximum and minimum singular values,  
      in relation to the mesh size, degree and regularity.
Moreover, in this paper we also study the effect of the number and 
      distribution of the collocation points on the spectral properties of 
      the collocation matrix, 
      providing insights into the optimization of the collocation points 
      for achieving better-conditioned linear systems.

\end{abstract}

\begin{keyword} 
	Isogeometric analysis \sep Isogeometric least-squares collocation method 
	\sep Condition number \sep Spectral properties
\end{keyword}

\end{frontmatter}


\section{Introduction}
\label{sec:intro}
Isogeometric analysis (IGA), introduced in \citep{Hughes2005IGA}, 
      is an emerging numerical technique providing connection between 
      computer-aided design (CAD) and computer-aided engineering (CAE).
The main idea of IGA is to directly handle the CAD geometry represented by  
      non-uniform rational B-splines (NURBS), and apply the same basis 
      functions for numerical simulation, 
      thereby achieving seamless integration of CAD and CAE.
Since Hughes et al. introduced IGA \citep{Hughes2005IGA}, 
      extensive studies have shown its advantages 
      in many problems and applications~\citep{bazilevs2006isogeometric,
      cottrell2009isogeometric,da2014mathematical}.
Compared to finite element method,
      IGA offers not only the standard $h-$ and $p-$refinement, 
      but also $k-$refinement, which enables smoother basis functions 
      and higher accuracy \citep{hughes2008duality}.

Recent studies have explored the isogeometric analysis collocation method (IGA-C) 
      to address the challenge of obtaining sparser stiffness matrices 
      compared to those arising from traditional 
      IGA Galerkin methods~\citep{Auricchio2011IGAC,
      auricchio2012isogeometric,manni2015isogeometric,anitescu}.
IGA-C involves substituting 
      collocation points into the strong form of a partial differential equation (PDE), 
      which results in a simple and flexible formulation 
      of the collocation matrix.
Moreover, since the collocation matrix requires only 
      one function evaluation per collocation point, 
      the computational cost of IGA-C is relatively low.
The number of collocation points can be set 
      equal to or greater than the degrees of freedom ($dof$),
      with the latter referred to as the isogeometric least-squares 
      collocation method (IGA-L)~\citep{anitescu,LinIGA-L}.

For collocation methods, the distribution of the collocation points 
      plays a crucial role in 
      improving the stability and convergence rate. 
An appropriate distribution of these points can enhance 
      the numerical performance significantly. 
Popular options, such as Greville abscissae and Demko
      ~\citep{Auricchio2011IGAC,demko}, 
      superconvergent (SC) points~\citep{anitescu},
      Cauchy-Galerkin (CG) points~\citep{Gomez}
      have been widely utilized and provided significant improvements 
      since their introduction.

Since the introduction of the IGA, 
      most of the studies have primarily focused on
      its applications \citep{bazilevs2006isogeometric,
      cottrell2009isogeometric,da2014mathematical,
      manni2011generalized,bazilevs2010isogeometric}, 
      with relatively little attention given to the theoretical estimation 
      of the spectral properties of the IGA matrices, such as 
      the condition number.
However, the accuracy and computational cost of solving the linear system 
      arising from the isogeometric discretization are crucial, 
      especially for large-scale problems,
      with the condition number playing a key role in both.
For direct solvers, an increase in the condition number 
      leads to a decrease in the number of significant digits in the solution, 
      resulting in a less accurate solution \citep{bathe2006finite}.
As the degree $p$ increases,
      the discretization matrix $A$ becomes denser,
      which increases the computational cost of direct solvers, 
      making them prohibitively expensive, 
      particularly for large problems~\citep{bathe2006finite}.
Consequently, iterative methods are typically employed for practical solutions.
For an iterative solver, the number of iterations 
      required to reach a specified error threshold 
      is also influenced by the condition number.
Systems with large condition numbers tend to require 
      more iterations to converge to the desired accuracy~\citep{bathe2006finite}.
Since the convergence rate of the iterative method 
      is strongly dependent on the condition number, 
      which ultimately impacts the accuracy of the numerical method, 
      it becomes essential to thoroughly investigate the relationship 
      between the condition number and the parameters of the IGA discretized model.

In this study, we aim to systematically investigate 
      the spectral properties of the collocation and mass matrix in IGA-L 
      for Poisson problem with homogeneous Dirichlet boundary conditions 
      in one-, two-, and three-dimensional domains.
Specifically, we analyze their behavior with respect to key parameters such as 
      the mesh size $h$, degree $p$, spatial dimension $d$,
      the number and the distribution of collocation points.

The rest of the paper is organized as follows.
Section \ref{sec:related} presents related studies. 
In Section \ref{sec:pre} we briefly review the fundamentals of B-splines and NURBS, 
      and introduce the eigenvalue and condition number estimations for the 
      IGA-Galerkin (IGA-G) approximation of the Poisson problem.
Section \ref{sec:IGAmodel} introduces the Poisson equation 
      with Dirichlet boundary conditions and its isogeometric discretization.
In Section \ref{sec:result}, we analyze the behavior of the 
      singular values and condition numbers of 
      the collocation and mass matrix, 
      providing estimations regarding their dependence on mesh size, degree 
      and regularity. 
      
\subsection{Related work}
\label{sec:related}
Despite the growing applications of IGA, 
      relatively few studies have explored 
      the theoretical properties of the IGA discretization matrix.
Gahalaut and Tomar derived bounds for the external eigenvalues 
      and the spectral condition number of the mass 
      and stiffness matrix in the isogeometric discretization 
      of the Poisson problem~\citep{gahalaut2012condition}, i.e.,
\begin{equation}
      \begin{aligned}
            &\mathcal{K}(\mathbf{M})\leq cp^{2}16^{p}\\
            &\mathcal{K}(\mathbf{K})\leq cp^{8}16^{p},
      \end{aligned}
\end{equation}
where $c$ is independent of $h$ and $p$, 
      and $\mathcal{K}(\mathbf{M})$, $\mathcal{K}(\mathbf{K})$ 
      denote the condition numbers of the mass matrix $\mathbf{M}$, 
      and stiffness matrix $\mathbf{K}$, respectively. 
Garoni et al. studied the spectral properties of the stiffness 
      matrix for classical second-order elliptic problems~\citep{garoni2014spectrum}.
They employed the concept of a symbol in the Toeplitz setting to
      examine the non-singularity, conditioning and 
      spectral distribution in the Weyl sense.
Eisentr$\ddot{a}$ger et al. investigated the impact of the choice of the shape functions 
      and element distortion on the condition number 
      of the stiffness matrix~\citep{eisentrager2020condition}.
Their study compared three high-order finite element methods (FEMs), i.e., 
      the $p$-version of the FEM \citep{szabo2021finite}, 
      the spectral element method (SEM) \citep{canuto2007spectral}, 
      and the NURBS-based IGA \citep{Hughes2005IGA}, 
      by comparing numerical results with the condition number estimations 
      from existing literature \citep{solin2003higher,demkowicz2006computing,melenk2002condition}.
Gervasio et al. conducted a systematic comparison of the theoretical properties 
      of SEM and NURBS-based IGA Galerkin methods 
      for approximating the Poisson problem~\citep{gervasio2020computational}.
Their analysis focused on the convergence properties, 
      algebraic structure, and spectral properties of the corresponding discrete systems. 
Additionally, the researchers proposed estimations for the behavior of the relevant 
      theoretical laws with respect to design parameters such as mesh size, 
      local degree, smoothness of NURBS basis functions, 
      spatial dimension, and the $dof$ involved in the computations.

For the IGA-C discretization matrices, 
      Zampieri et al. analyzed the spectral properties of the mass and 
      collocation matrices for specific acoustic wave problems,  
      where they employed the IGA-C in space and 
      the Newmark method in time~\citep{zampieri2024conditioning}. 
However, a comprehensive study and estimation 
      of the spectral properties of the IGA-L discretization 
      matrices for standard PDEs is still lacking.

\section{Preliminaries}
\label{sec:pre}
  \subsection{B-spline and NURBS basis functions}
Consider a set of distinct non-decreasing real numbers 
      $\left\{\xi_{0}=0,\xi_{1},\cdots,\xi_{n-1},\xi_{n}=1\right\}$
      in the one-dimensional patch $\left[0,1\right]$,
      and given two positive integers $p$ and $k$
      with $0\leq k \leq p-1$,
      let
      \begin{equation}
      \label{eq:knot}
            \{\underbrace{\xi_{0},\cdots,\xi_{0}}_{p+1},
            \underbrace{\xi_{1},\cdots,\xi_{1}}_{p-k},\cdots,
            \underbrace{\xi_{n-1},\cdots,\xi_{n-1}}_{p-k},
            \underbrace{\xi_{n},\cdots,\xi_{n}}_{p+1}\}
      \end{equation}
      be an open knot vector, i.e., the boundary knots have 
      $p+1$ multiplicity.
The $i-$th univariate B-spline basis function $B_{i}^{p}(\xi)$ with the 
      degree $p\geq 1$ and the regularity $C^{k}$
      can be defined over the above knot vector by means of 
      the deBoor-Cox recursion formula~\citep{piegl2012nurbs}. 
The number of linear dependent basis functions is $N_{b}=(n-1)(p-k)+p+1$.
In this paper, we utilize the uniform open knot vector, 
      in which the internal knot values are uniformly distributed.
Therefore, the mesh size $h=\frac{1}{n}$.

The regularity is a prominent property of B-splines, 
      which denotes the global continuity of the basis functions.
If internal knots are repeated $p-k$ times,
      then the B-splines are $C^{k}$-continuous,
      that stands for the \textit{global order of regularity} \citep{cottrell2009isogeometric}.
In this paper, we consider two extreme values for $k$. 
When $k=p-1$, the B-splines are globally $C^{p-1}$, and 
      we use the notation $IGA-C^{p-1}$ to identify the case.
When $k=1$, we use the notation $IGA-C^{1}$.

In $d-$dimensional cases, the B-spline basis functions 
      are calculated by tensor product.
We introduce the parametric domain $\hat{\Omega}=[0,1]^{d}$ 
      with a knot vector (\ref{eq:knot}) 
      in each parametric direction.
The multivariate B-spline basis on $\hat{\Omega}$ is 
      $(B_{i_{1}}^{p_{1}}\cdots B_{i_{d}}^{p_{d}})$.
      
NURBS basis functions are defined starting from B-spline basis functions 
      by associating a set of weights 
      $\{\omega_{1},\omega_{2},\cdots,\omega_{N_{b}}\}$.
In this paper, we assume $\omega_{i}\in \mathbb{R},\omega_{i}>0$ 
      for $i=1,2,\cdots,N_{b}$.
The $i$-th univariate NURBS basis function is
\begin{equation}
      \label{eq:nurbs}
      N_{i}^{p}(\xi)=\frac{\omega_{i}B_{i}^{p}(\xi)}
            {\sum_{j=1}^{N_{b}}\omega_{j}B_{j}^{p}(\xi)}.
\end{equation}

NURBS inherit properties from their B-splines counterpart, 
      so the degree of corresponding NURBS is $p$,
      and the global regularity is $k$.

\subsection{Spectral properties estimations of IGA-G}
The eigenvalues and the condition number of IGA-G 
      matrices in Poisson equation have been systematically studied 
      in~\citep{gervasio2020computational}.
Let $\mathbf{M}=(m_{ij})$ be the mass matrix, where 

\begin{equation}
      \label{eq:mass_galerkin}
      m_{ij}=(N_{i}^{p},N_{j}^{p})=\int_{\Omega}N_{i}^{p}N_{j}^{p}.
\end{equation} 
Let $\mathbf{K}=(k_{ij})$ be the stiffness matrix, where

\begin{equation}
      \label{eq:stiff_galerkin}
      k_{ij}=\int_{\Omega}\nabla N_{i}^{p}\cdot \nabla N_{j}^{p}.
\end{equation} 

The estimations in ~\citep{gervasio2020computational} show that 
      for $k=0$ regularity, it holds:
      \begin{equation}
            \label{eq:estimate_IGA-G_0}
            \begin{aligned}
                  &\lambda_{min}(\mathbf{M}_{0})\thicksim h^{d}p^{-d/2}4^{-pd}\\
                  &\lambda_{max}(\mathbf{M}_{0})\thicksim h^{d}p^{-d}\\
                  &\mathcal{K}(\mathbf{M}_{0})\thicksim p^{-d/2}4^{pd}\\
                  &\lambda_{min}(\mathbf{K}_{0})\thicksim  
                        \begin{cases}
                              h^{d}p^{-d} &\text{if} \ h 
                                    \lesssim (p^{2+d/2}4^{-dp})^{1/2}\\
                              h^{d-2}p^{2-d/2}4^{-dp} &\text{otherwise}
                        \end{cases}\\
                  &\lambda_{max}(\mathbf{K}_{0})\thicksim h^{d-2}p^{2-d}\\
                  &\mathcal{K}(\mathbf{K}_{0})\thicksim 
                        \begin{cases}
                              h^{-2}p^{2} &\text{if} \ h 
                                    \lesssim (p^{2+d/2}4^{-dp})^{1/2}\\
                              p^{-d/2}4^{dp} &\text{otherwise}
                        \end{cases}\\
            \end{aligned}
      \end{equation}
      for $d=1,2,3$, respectively. 
$\mathbf{M}_{0}$ represent the mass matrix  
      and $\mathbf{K}_{0}$ is the stiffness matrix for $k=0$ regularity.
Let $\lambda_{min}(M_{0})$, $\lambda_{max}(M_{0})$ represent the 
      minimum and maximum eigenvalues of $\mathbf{M}_{0}$, respectively.
The condition number of $\mathbf{M}_{0}$ is denoted by $\mathcal{K}(M_{0})$.
The meanings of $\lambda_{min}(\mathbf{K}_{0})$, $\lambda_{max}(\mathbf{K}_{0})$ 
      and $\mathcal{K}(\mathbf{K}_{0})$ are analogous to 
      those corresponding notations for $\mathbf{M}_{0}$.

We denote with $\mathbf{M}_{p-1}$ the mass matrix 
      and $\mathbf{K}_{p-1}$ the stiffness matrix
      associated with IGA-$C^{p-1}$ approximation.
For $k=p-1$ regularity, it holds:
      \begin{align}
            \label{eq:estimate_IGA-G_p}
                  &\lambda_{min}(\mathbf{M}_{p-1})\thicksim 
                        \begin{cases}
                              h^{d}e^{-pd} &\text{if} \ h \lesssim 1/p\\
                              (\frac{e}{4})^{-d/h}(\frac{h}{p})^{d/2}4^{-pd} 
                                    &\text{otherwise}
                        \end{cases} \notag\\
                  &\lambda_{max}(\mathbf{M}_{p-1})\thicksim 
                        \begin{cases}
                              h^{d} &\text{if} \ h \lesssim 1/p\\
                              p^{-d} &\text{otherwise}
                        \end{cases} \notag\\
                  &\mathcal{K}(\mathbf{M}_{p-1})\thicksim 
                        \begin{cases}
                              e^{pd} &\text{if} \ h \lesssim 1/p\\
                              (\frac{e}{4})^{d/h}(hp)^{-d/2}4^{pd} &\text{otherwise}
                        \end{cases} \notag \displaybreak\\
                  &\lambda_{min}(\mathbf{K}_{p-1})\thicksim  
                        \begin{cases}
                              h^{d} &\text{if} \ h \lesssim e^{-pd/2}\\
                              h^{d-2}e^{-pd} &\text{if} \ e^{-pd/2} 
                                    \lesssim h \lesssim 1/p\\
                              (\frac{e}{4})^{-d/h}p^{2-d/2}h^{d/2}4^{-pd} 
                                    &if \ h \gtrsim 1/p
                        \end{cases}\\
                  &\lambda_{max}(\mathbf{K}_{p-1})\thicksim 
                        \begin{cases}
                              ph^{d-2} &\text{if} \ h \lesssim 1/p \quad for \ p>2\\
                              p^{2-d}h^{-1} &\text{otherwise}
                        \end{cases}  \notag\\
                  &\mathcal{K}(\mathbf{K}_{p-1})\thicksim 
                        \begin{cases}
                              h^{-2}p &\text{if} \ h \lesssim e^{-pd/2}\\
                              pe^{dp} &\text{if} \ e^{-pd/2} 
                                    \lesssim h \lesssim 1/p\\
                              (\frac{e}{4})^{d/h}p^{-d/2}h^{-d/2-1}4^{dp} 
                                    &\text{otherwise}.
                        \end{cases}  \notag\\
                  \nonumber
      \end{align}
In this paper, we study the spectral properties 
      of the mass and collocation matrices arising in IGA-L 
      for Poisson problems, 
      and compare the numerical results with 
      the estimations of IGA-G in the literature 
      \citep{gervasio2020computational}.
      
\section{IGA-L discretization model}
\label{sec:IGAmodel}
Let $\Omega \subset \mathbb{R}^{d}$, with $d=1,2,3$, be a bounded domain,
      and let $f$ be a given function.
We consider the Poisson equation with homogeneous Dirichlet boundary conditions, 
\begin{equation}
	\label{eq:pde}
	\left\{
	\begin{array}{ll}
		-\Delta u=f & \ \text{in }\Omega \\
		u=0 & \ \text{on }\partial\Omega ,\\ 
	\end{array}
	\right.
\end{equation}
where $u$ is the unknown solution. 
We discretize the strong form of the Poisson problem (\ref{eq:pde})
      using the least-squares collocation method of the
      isogeometric analysis with NURBS basis functions.

\subsection{Isogeometric least-squares collocation method}      
In IGA, NURBS are used to represent both the computational domain 
      $\Omega \subset \mathbb{R}^{d}$ and the numerical solution.
The geometric mapping is 
      derived by associating each basis function with 
      a corresponding control point $\bm{P}_{i}\in \mathbb{R}^{d}$.
Therefore, every point $\bm{x}$ in the physical domain 
      $\Omega$ can be expressed through the geometric mapping $G$, 
\begin{equation}
     G:\bm{x}(\bm{\xi})=
      \sum_{i=1}^{N_{b}}\bm{P}_{i}N_{i}^{p}(\bm{\xi}).
      \label{eq:geometric_map}
\end{equation}

Assuming the geometric mapping is invertible, 
      the mesh in the physical domain $\Omega$ 
      can be uniquely determined by its corresponding
      mesh in the parameter domain $\hat{\Omega}$.
The space of NURBS functions in the physical domain 
      is defined as:
\begin{equation}
	\mathcal{V}_{N_{b}}=\text{span}\left\{N_{i}^{p}\circ G^{-1},
      i=1,\cdots,N_{b}\right\}.
      \nonumber
\end{equation}
In two or higher dimensions, the space can be defined analogically:
\begin{equation}
	\begin{aligned}
		\mathcal{V}_{N_{i},N_{j}}&=\text{span}\{N_{i,j}^{p,q}\circ G^{-1},
                  i=1,\cdots,N_{i};j=1,\cdots,N_{j}\} \\
		\mathcal{V}_{N_{i},N_{j},N_{k}}&=\text{span}\{N_{i,j,k}^{p,q,r}\circ G^{-1},
                  i=1,\cdots,N_{i};j=1,\cdots,N_{j};\\
			&\qquad\qquad\qquad\qquad\ k=1,\cdots,N_{k}\}.
            \nonumber
	\end{aligned}
\end{equation}

IGA seeks the numerical solution to the problem (\ref{eq:pde}) 
      within the NURBS space incorporating the boundary conditions 
      in the physical domain. 
The key idea behind IGA-L is to sample two sets of collocation points:
      $\left\{\tau_{j_{1}}^{in},j_{1}=1,\cdots,m^{in}\right\}$ 
      in the interior of the physical domain $\Omega$, 
      and $\left\{\tau_{j_{2}}^{bd},j_{2}=1,\cdots,m^{bd}\right\}$ 
      on the boundary $\partial\Omega$.  
Thus, solving the IGA-L approximation of (\ref{eq:pde}) is equivalent to
      finding the best approximation 
      $U_{h}\in \mathcal{V}_{N_{b}}$ such that
	\begin{equation}
	      \label{eq:discretized}
	      \left\{
	      \begin{array}{ll}
		      -\Delta U_{h}(\tau_{j_{1}}^{in})=f(\tau_{j_{1}}^{in})
                  & \  j_{1}=1,\cdots ,m^{in} \vspace{0.5em} \\
		      U_{h}(\tau_{j_{2}}^{bd})=0
                  & \  j_{2}=1,\cdots,m^{bd}. \\  
	      \end{array}
	      \right.
      \end{equation}

The collocation points are determined based on the knot vectors.
A commonly used set of collocation points is Greville abscissae~\citep{Auricchio2011IGAC}. 
Given a knot vector with regularity $k=p-1$ as defined in (\ref{eq:knot}), i.e., 
\begin{equation}
      \nonumber
      \{\underbrace{\xi_{0},\cdots,\xi_{0}}_{p+1},
            \xi_{1},\cdots,
            \xi_{n-1},
            \underbrace{\xi_{n},\cdots,\xi_{n}}_{p+1}\}, 
\end{equation}
            we denote it as 
            $\{\eta_{0},\eta_{1},\cdots,\eta_{n+2p}\}$ for ease of presentation. 
The Greville abscissae points $\hat{\xi_{i}}, i=0,1,\cdots,n+p-1$ for classical 
      IGA-C in the parameter domain are obtained by:
\begin{equation}
	\hat{\xi_{i}}=\frac{\eta_{i+1}+\cdots+\eta_{i+p}}{p},\quad i=0,\cdots,n+p-1.
\end{equation}
The corresponding collocation points in the physical domain 
      are obtained through the geometric mapping (\ref{eq:geometric_map}): $\tau_{i}=x(\hat{\xi_{i}})$. 
Let $m$ denotes the total number of collocation points in IGA-L. 
To generate these $m$ points, 
      a knot vector of size $m+p+1$ with regularity $k=p-1$ 
      is constructed, i.e., 
\begin{equation}
      \nonumber
      \{\underbrace{\xi_{0},\cdots,\xi_{0}}_{p+1},
            \xi_{1},\cdots,
            \xi_{m-p-1},
            \underbrace{\xi_{m-p},\cdots,\xi_{m-p}}_{p+1}\}, 
\end{equation}
      which can be denoted as 
      $\{\eta_{0},\eta_{1},\cdots,\eta_{m+p}\}$.
The Greville abscissae points are then calculated as follows:
\begin{equation}
	\hat{\xi_{i}}=\frac{\eta_{i+1}+\cdots+\eta_{i+p}}{p},\quad i=0,\cdots,m-1.
\end{equation}
The collocation points in the physical domain are given by 
      the geometric mapping (\ref{eq:geometric_map}): $\tau_{i}=x(\hat{\xi_{i}})$.

SC points are developed based on the superconvergent theory.
Anitescu et al. aimed to identify collocation points 
      where the numerical solution exhibits behavior 
      similar to that of the standard Galerkin solution~\citep{anitescu}.
Table \ref{tbl:scpoints} presents the locations of superconvergent points within  
the reference interval $\left[-1,1\right]$ for  
different degrees.
These SC points are obtained by mapping the superconvergent points 
      from the reference interval to 
      the knot interval $\left[\xi_{i},\xi_{i+1}\right]$.
Since the number of superconvergent points is approximately 
      twice the number of $dof$,
      and all these points are utilized,
      IGA-C with SC points naturally leads to the formation 
      of a least-squares system.

      \begin{table}[h!]
            \renewcommand{\arraystretch}{1.2}
            \captionsetup{singlelinecheck=false,width=0.88\linewidth,
                          justification=raggedright}
            \caption{Superconvergent points in the reference interval.}
            \label{tbl:scpoints}
            \centering
            \begin{threeparttable}
                  \begin{tabular}{ccccc}
                        \hline
                        &\textbf{Degree} & &\textbf{Superconvergent Points}&\\
                        \hline               
                        &$p=3$  & &$-\frac{1}{\sqrt{3}},\frac{1}{\sqrt{3}}$&\\
                        &$p=4$  & &$-1,0,1$&\\
                        &$p=5$  & &$-\frac{\sqrt{225-30\sqrt{3}}}{15},
                                    \frac{\sqrt{225-30\sqrt{3}}}{15}$&\\
                        &$p=6$  & &$-1,0,1$&\\
                        &$p=7$  & &$-0.50491856751,0.50491856751$&\\
                        \hline                                             
                  \end{tabular}   
            \end{threeparttable}       
      \end{table}

\begin{figure}[htbp]
      \centering
      \subfigure[]{
            \label{subfig:dis_gre_p3}
            \includegraphics[width=0.92\linewidth]{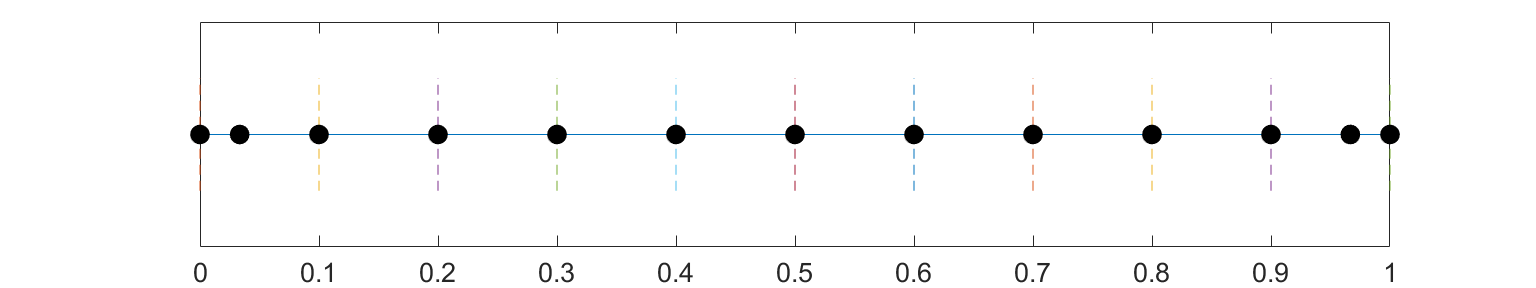}}\\
      \subfigure[]{
            \label{subfig:dis_SC_p3}
            \includegraphics[width=0.92\linewidth]{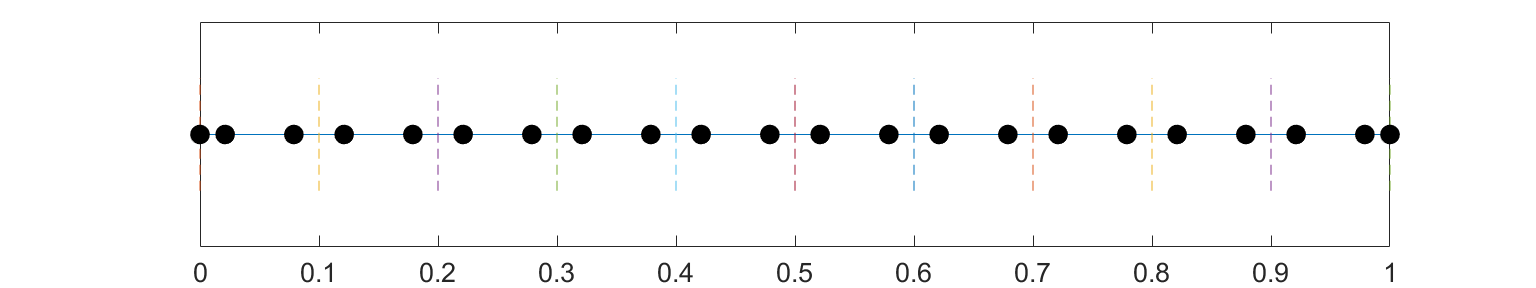}}\\
      \subfigure[]{
            \label{subfig:dis_CG_p3}
            \includegraphics[width=0.92\linewidth]{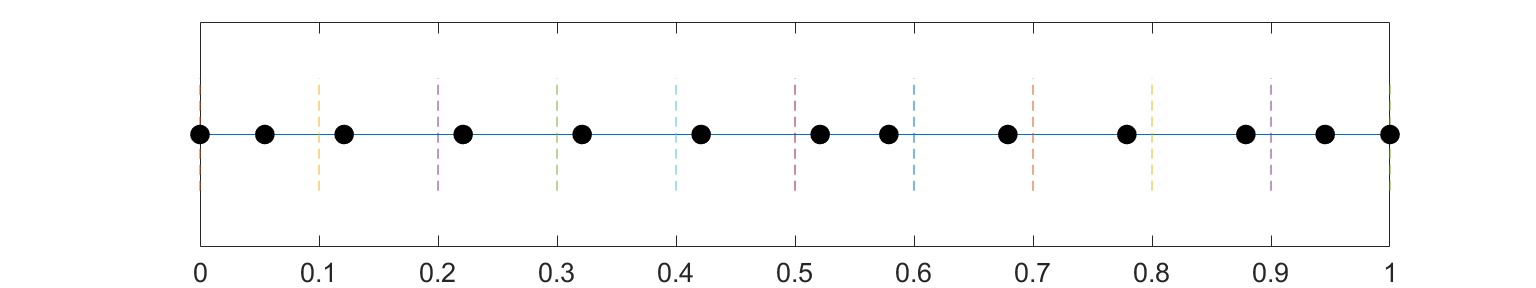}}\\
      \caption{Distributions of different types of collocation points ($p=3$). 
            $(\mathbf{a})$ Greville points.  $(\mathbf{b})$ SC points. 
            $(\mathbf{c})$ CG points.}
      \label{fig:dis_collocation}
\end{figure}

CG points are selected as a specific subset of the 
      superconvergent points~\citep{Gomez}.
They retain one superconvergent point in each knot interval 
      while maintaining global symmetry.
As a result, the number of CG points equals the number of $dof$.
Fig.\ref{fig:dis_collocation} shows the distribution 
      of Greville points, SC points and CG points 
      when the degree $p=3$.

\subsection{Matrix form of IGA-L}
In IGA, the numerical solution $U_{h}$ is defined on a mesh 
      determined by the geometric mapping.
Consequently, $U_{h}$ in IGA is represented as a linear combination 
      of the NURBS basis functions in $\Omega$, 
      which also serve as the basis functions for the geometric mapping $G$.
      
For one-dimensional problems, 
      the numerical solution can be expressed as:
\begin{equation}
	\label{eq:num_sol}
	U_{h}(\tau)=\sum_{i=1}^{N_{b}}u_{i}N_{i,p}(G^{-1}(\tau)),
      \quad \tau\in\Omega
\end{equation}
where $\{u_{i}\}_{i=1}^{N_{b}}$ are the unknown control coefficients. 
By substituting (\ref{eq:num_sol}) to (\ref{eq:discretized}) 
      we obtain a non-symmetric linear system 
      with respect to the unknowns $\{u_{i}\}_{i=1}^{N_{b}}$. 
The vector $\{u_{i}\}_{i=1}^{N_{b}}$ is arranged into a column vector, 
      and the basis functions and collocation points are 
      reordered such that the first $N_{b}^{in}$ basis functions
      correspond to the internal control coefficients, 
      and the first $m^{in}$ collocation points 
      lie within the support of the internal basis functions.

We define $\mathbf{u}=\left[u_{i}\right]_{i=1}^{N_{b}^{in}}$,      
      and $\mathbf{b}=\left[f(\tau_{j})\right]_{j=1}^{m^{in}}$.
Thus, the problem (\ref{eq:discretized}) can be 
      reformulated in a matrix form as:
\begin{equation}
	\label{eq:linear_equation}
	\mathbf{Au}=\mathbf{b}
\end{equation}
with 
\begin{align}
	\label{eq:IGAC_model}	
      \nonumber
	&\mathbf{A}_{ji}=-\Delta N_{i,p}(G^{-1}(\tau_{j})) \quad j=1,\cdots,m^{in};
            \ i=1,\cdots,N_{b}^{in},\\
      &\mathbf{b}_{j}=f(\tau_{j}) \quad j=1,\cdots,m^{in}.	
      \nonumber
\end{align}

	
$\mathbf{A}$ denotes the \emph{collocation matrix} with $m^{in}$ rows 
      and $N_{b}^{in}$ columns, 
      and $\mathbf{b}$ is the \emph{load vector}. 
The \emph{mass matrix} $\mathbf{M}$ in IGA-L is defined as 
a collocation matrix on the basis functions directly: 
\begin{equation}
	\label{eq:mass matrix}
	\mathbf{M}_{ji}=N_{i,p}(G^{-1}(\tau_{j})) \quad j=1,\cdots,m^{in};
      \ i=1,\cdots,N_{b}^{in}.
\end{equation}

The numerical solution is then obtained by 
      solving the least-squares linear system (\ref{eq:linear_equation}).
The corresponding least-squares problem is given by:
\begin{equation}
	\label{eq:ls_problem}
	\mathop{\min}_{\mathbf{u}}\|\mathbf{Au}-\mathbf{b}\|^{2}.
\end{equation}
The normal equations associated with (\ref{eq:ls_problem}) are:
\begin{equation}
	\label{eq:normal_equation}
	\mathbf{A}^{\text{T}}\mathbf{A}\mathbf{u}
      =\mathbf{A}^{\text{T}} \mathbf{b}.
\end{equation}

\section{Numerical analysis}
\label{sec:result}
In this section, we analyze the condition number and singular values 
      of the collocation and mass matrices in IGA-L model 
      for the Poisson problem with homogeneous 
      Dirichlet boundary conditions.

For any matrix $\mathbf{A}$ (which may not be square), 
      let $\sigma_{max}(\mathbf{A})$ and $\sigma_{min}(\mathbf{A})$ denote 
      its maximum and minimum singular values, respectively.
The $spectral\ condition\ number$ of $\mathbf{A}$ is given by:
      \begin{equation}
            \label{eq:condition number}
            \mathcal{K}(\mathbf{A})
            =\frac{\sigma_{max}(\mathbf{A})}{\sigma_{min}(\mathbf{A})}.
      \end{equation} 
The spectral condition number of the matrix $\mathbf{A}$ increases as  
      its columns or rows become more unevenly distributed or 
      closer to linear dependence. 
The former tends to increase $\sigma_{max}(\mathbf{A})$, 
      while the latter tends to decrease $\sigma_{min}(\mathbf{A})$.

Based on Section \ref{sec:IGAmodel}, we summarize the 
      factors affecting the spectral properties of the 
      collocation and mass matrices in IGA-L as follows:
      \begin{itemize}
            \item The type of partial differential equations, 
                  including the differential operators and 
                  the boundary conditions.
            \item Geometric mapping, i.e., the parameterization of 
                  the physical domain.
            \item Dimension $d$.
            \item Degree $p$ of the NURBS basis functions.
            \item Mesh size $h$.
            \item Regularity $k$.
            \item Number of collocation points $m$.
            \item Distribution of collocation points.
      \end{itemize}
For the numerical analysis, 
      we primarily focus on the effects of the degree $p$, 
      mesh size $h$, regularity $k$, dimension $d$,
      as well as the number and distribution of collocation points 
      on the spectral properties of the collocation and mass matrices.
We present estimates for the condition number 
      and singular values of IGA-L matrices with respect to $h$ and $p$, 
      and explore how the number and distribution of the collocation points 
      influence these properties.
      
\begin{figure}[H]
	\centering
	\subfigure[A quarter of the annulus]{
            \label{subfig:2d_domain}
		\includegraphics[width=0.43\linewidth]{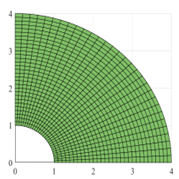}
	}
	\quad
	\subfigure[A unit cube]{
            \label{subfig:cube_domain}
		\includegraphics[width=0.48\linewidth]{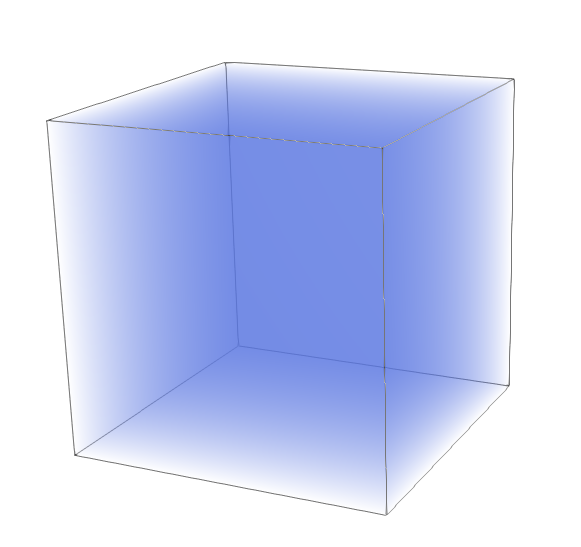}
	}
	\quad
	\centering
	\subfigure[A $1/8$th of a hollow sphere]{
            \label{subfig:hollow_sphere}
		\includegraphics[width=0.52\linewidth]{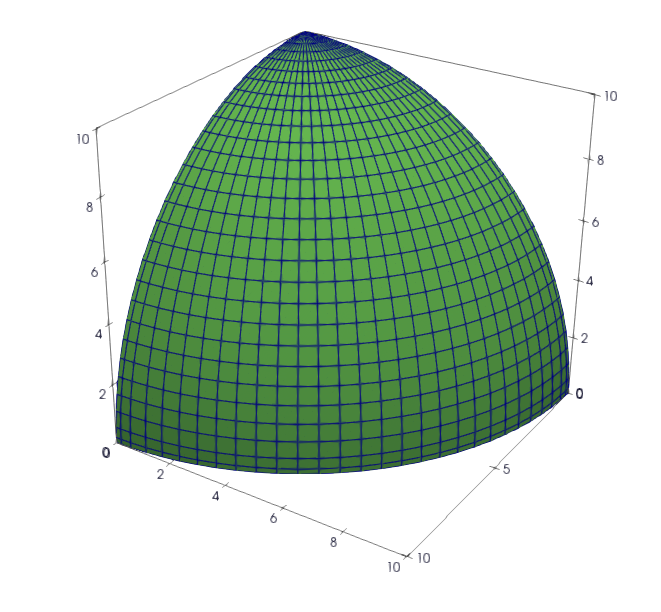}
	}
	\quad
	\caption{Physical domain in two- and three-dimensional cases.}
      \label{fig:phy_domain}
\end{figure}

We consider the Poisson problem 
      in one-, two-, and three-dimensional physical domains.
The one-dimensional physical domain is $\Omega=\left[0,1\right]$.
The two-dimensional physical domain is a quarter of the annulus,  
      and for the three-dimensional case, we study 
      a unit cube $\Omega=\left[0,1\right]^{3}$ 
      and a model of $1/8$th of a hollow sphere.
The radius of the mid-surface of the hollow sphere is $10$ 
      and the thickness of the model is $0.04$.
Fig. \ref{fig:phy_domain} shows the two- and three-dimensional 
      physical domains we study in this paper.

For the regularity, we consider two cases: $k=p-1$ (denoted as $IGA-C^{p-1}$) 
      and $k=1$ (denoted as $IGA-C^{1}$). 
      The corresponding collocation matrices are denoted as 
      $\mathbf{A}_{p-1}$ and $\mathbf{A}_{1}$, respectively.
Similarly, the mass matrices for these two cases are 
      denoted as $\mathbf{M}_{p-1}$ and $\mathbf{M}_{1}$, respectively.

Due to RAM and memory limitations, 
      the ranges of $p$ and $h$ vary 
      depending on the regularity and the spatial dimension.
For the distribution of collocation points, 
      we consider the Greville points, 
      SC points and CG points.
We first focus on the Greville points in the following analysis, 
      where the number of collocation points is 
      initially set to four times the $dof$.
We summarize the main results, conduct numerical analysis of both the 
      collocation matrix and the mass matrix in IGA-L, 
      and present estimations for their 
      condition number and singular values in relation to $p$, $h$ and $k$.
All the experiments are conducted on a PC  
      with an Intel Core i7-10700 2.90GHz CPU and 16GB of RAM.

\subsection{IGA-$C^{p-1}$ collocation matrix} 
In this section, we study the spectral properties of the IGA-$C^{p-1}$ 
      collocation matrix $\mathbf{A}_{p-1}$.
The mesh size $h$ ranges from $0.1$ to $0.01$, 
      with a decrease from $0.2$ to $0.05$ when $d=3$.
The degree $p$ starts at $2$, 
      and increases to $24,16,7$ for $d=1,2,3$, respectively. 
For the hollow sphere model, the mesh size ranges from $0.1$ to $0.05$, 
      and the degree $p$ starts at $4$ and increases to $7$.

The numerical results show that for $h>0$, $p\geq 2$ and $d=1,2,3$, 
      $\sigma_{min}(\mathbf{A}_{p-1})$, $\sigma_{max}(\mathbf{A}_{p-1})$ and 
      $\mathcal{K}(\mathbf{A}_{p-1})$ behave as:
\begin{align}       
      \label{eq:estimate_pA_min}
            &\sigma_{min}(\mathbf{A}_{p-1})\thicksim  
                  \begin{cases}
                        c_{1} &\text{if} \ h \lesssim e^{-pd/4}p^{d/8}\\
                        (e/2)^{-2dp}p^{(4/e)^{d}}h^{-2} 
                              &\text{ otherwise },     
                  \end{cases}\\  
            \label{eq:estimate_pA_max}
            &\sigma_{max}(\mathbf{A}_{p-1})\thicksim h^{-2}p^{2},\\
            \label{eq:estimate_pA_cond}
            &\mathcal{K}(\mathbf{A}_{p-1})\thicksim 
                  \begin{cases}
                        h^{-2}p^{2} &\text{if} \ h \lesssim e^{-pd/4}p^{d/8}\\
                        (e/2)^{2dp}p^{2-(4/e)^{d}} &\text{otherwise},
                  \end{cases}
\end{align}
where $c_{1}$ is independent of $h$ and $p$. 
The symbol $\thicksim$ means "\textit{up to a constant independent of both h and p}".

\begin{figure}[htbp]
		\centering 
		\subfigure[Condition number vs. $h$]{
			  \label{subfig:1d_p_condA_h}
			  \includegraphics[width=0.39\linewidth]{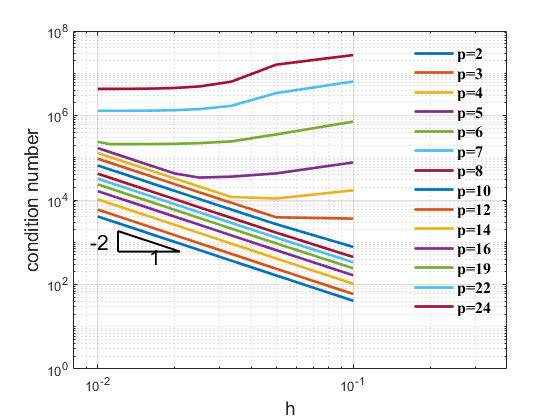}}    
		\subfigure[Condition number vs. $p$]{
			  \label{subfig:1d_p_condA_p}
			  \includegraphics[width=0.39\linewidth]{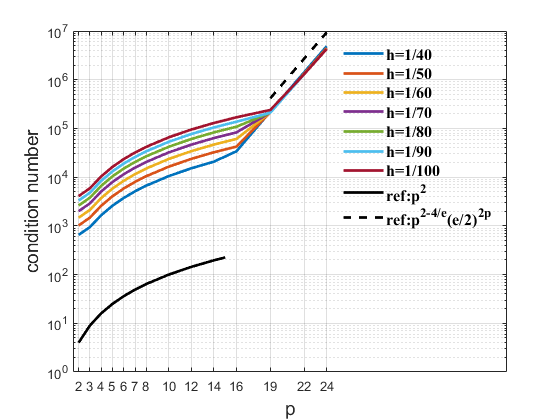}}\\             
		\subfigure[Maximum singular value vs. $h$]{
			  \label{subfig:1d_p_maxA_h}
			  \includegraphics[width=0.39\linewidth]{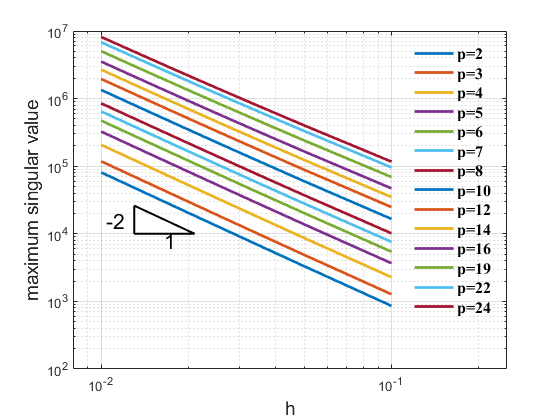}}                
		\subfigure[Maximum singular value vs. $p$]{
			  \label{subfig:1d_p_maxA_p}
			  \includegraphics[width=0.39\linewidth]{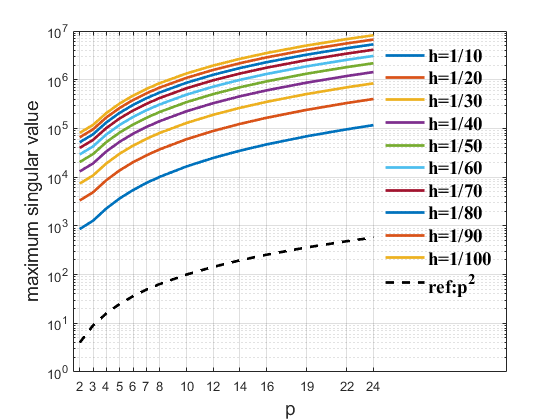}}\\
		\subfigure[Minimum singular value vs. $h$]{
			  \label{subfig:1d_p_minA_h}
			  \includegraphics[width=0.39\linewidth]{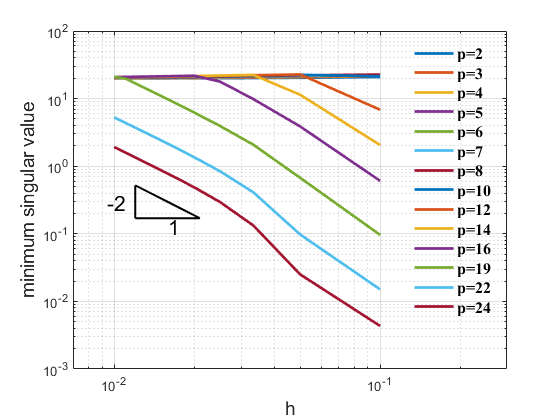}}               
		\subfigure[Minimum singular value vs. $p$]{
			  \label{subfig:1d_p_minA_p}
			  \includegraphics[width=0.39\linewidth]{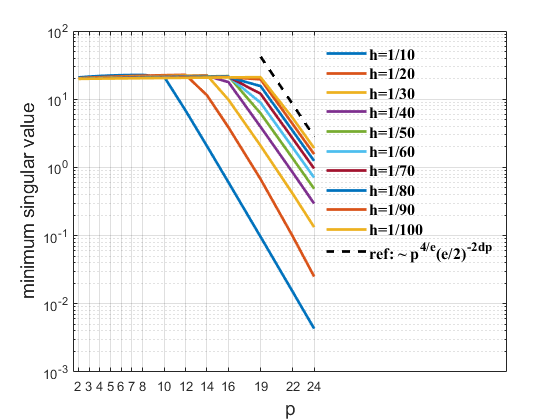}}
		\caption{The maximum/minimum singular value and the spectral condition number 
		of $\mathcal{K}(\mathbf{A}_{p-1})$ for $d=1$, versus $h$ (at left) and 
		versus $p$ (at right).}
		\label{fig:1d_p_collocation}
\end{figure}

\begin{figure}[htbp]
            \centering  
            \subfigure[Condition number vs. $h$]{
                 \label{subfig:2d_p_condA_h}
                  \includegraphics[width=0.39\linewidth]{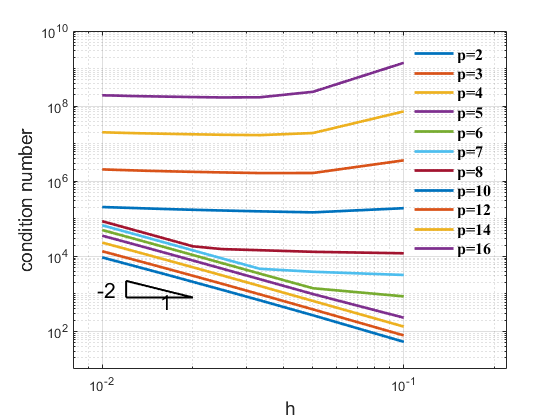}}    
            \subfigure[Condition number vs. $p$]{
                  \label{subfig:2d_p_condA_p}
                  \includegraphics[width=0.39\linewidth]{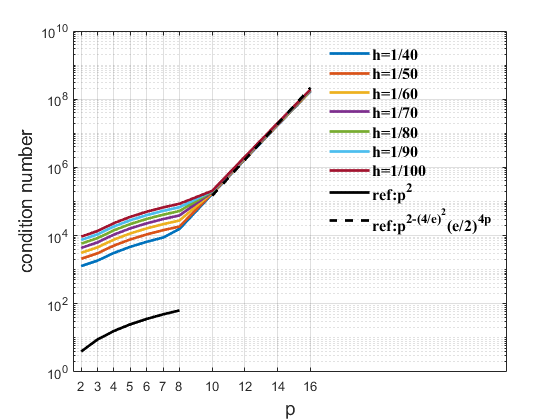}}\\             
            \subfigure[Maximum singular value vs. $h$]{
                  \label{subfig:2d_p_maxA_h}
                  \includegraphics[width=0.39\linewidth]{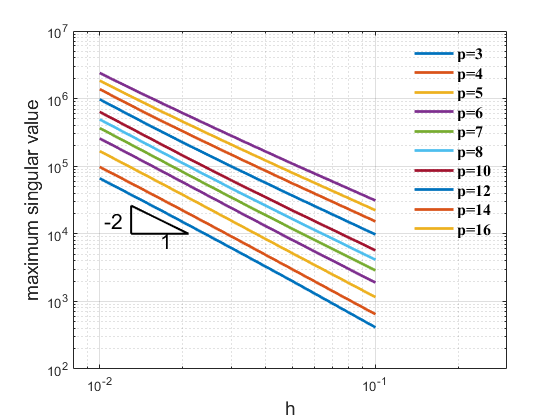}}                
            \subfigure[Maximum singular value vs. $p$]{
                  \label{subfig:2d_p_maxA_p}
                  \includegraphics[width=0.39\linewidth]{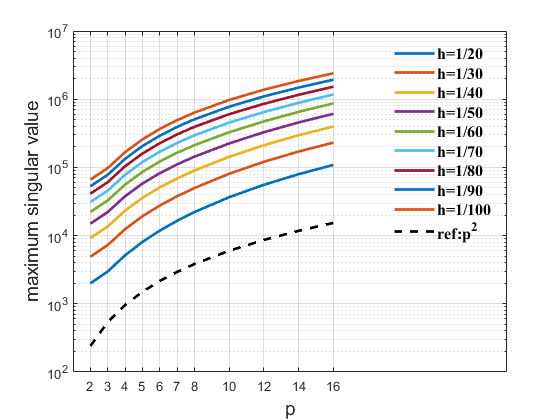}}\\
            \subfigure[Minimum singular value vs. $h$]{
                  \label{subfig:2d_p_minA_h}
                  \includegraphics[width=0.39\linewidth]{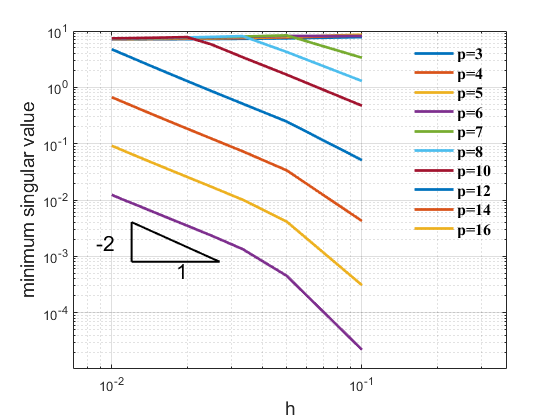}}               
            \subfigure[Minimum singular value vs. $p$]{
                  \label{subfig:2d_p_minA_p}
                  \includegraphics[width=0.39\linewidth]{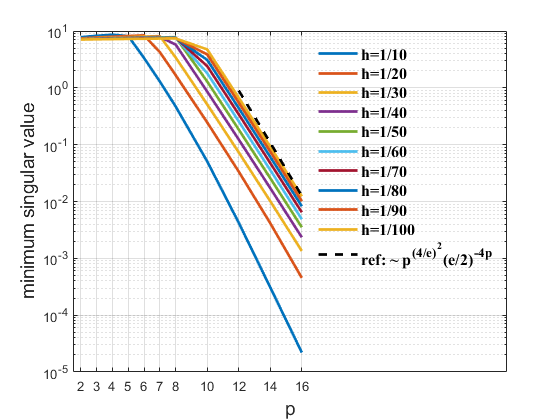}}
            \caption{The maximum/minimum singular value and the spectral condition number 
            of $\mathcal{K}(\mathbf{A}_{p-1})$ for $d=2$, versus $h$ (at left) and 
            versus $p$ (at right).}
            \label{fig:2d_p_collocation}
      \end{figure}

\begin{figure}[htbp]
      \centering  
      \subfigure[Condition number vs. $h$]{
            \label{subfig:3d_p_condA_h}
            \includegraphics[width=0.39\linewidth]{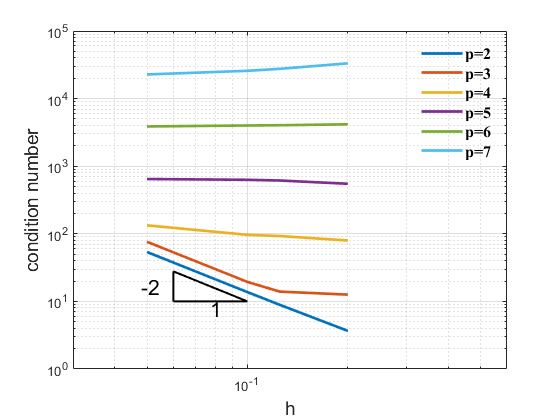}}    
      \subfigure[Condition number vs. $p$]{
            \label{subfig:3d_p_condA_p}
            \includegraphics[width=0.39\linewidth]{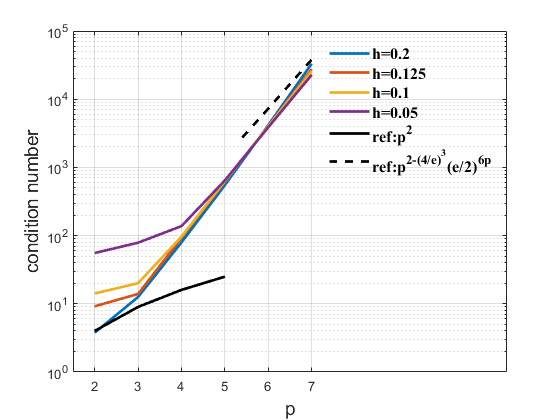}}\\             
      \subfigure[Maximum singular value vs. $h$]{
            \label{subfig:3d_p_maxA_h}
            \includegraphics[width=0.39\linewidth]{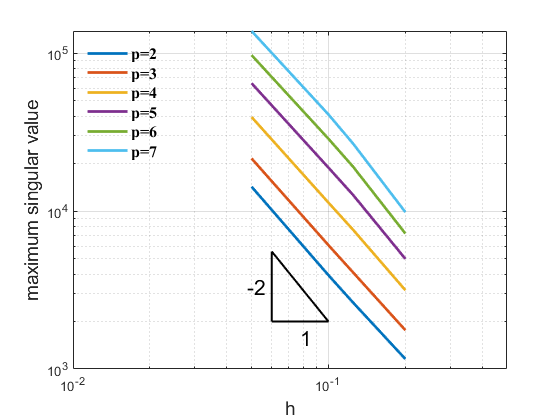}}                
      \subfigure[Maximum singular value vs. $p$]{
            \label{subfig:3d_p_maxA_p}
            \includegraphics[width=0.39\linewidth]{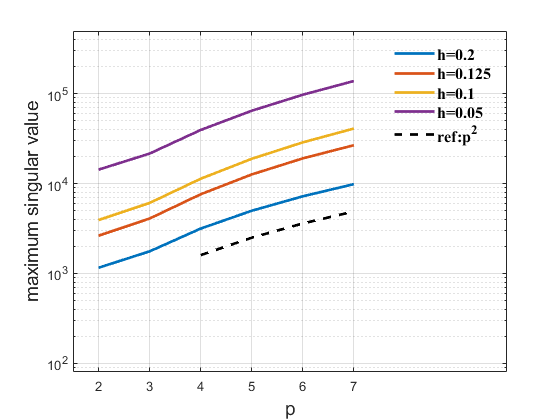}}\\
      \subfigure[Minimum singular value vs. $h$]{
            \label{subfig:3d_p_minA_h}
            \includegraphics[width=0.39\linewidth]{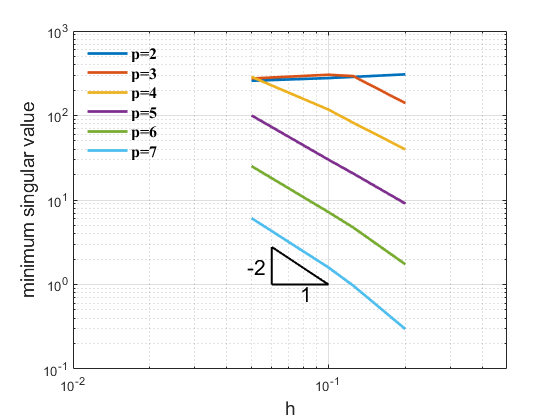}}               
      \subfigure[Minimum singular value vs. $p$]{
            \label{subfig:3d_p_minA_p}
            \includegraphics[width=0.39\linewidth]{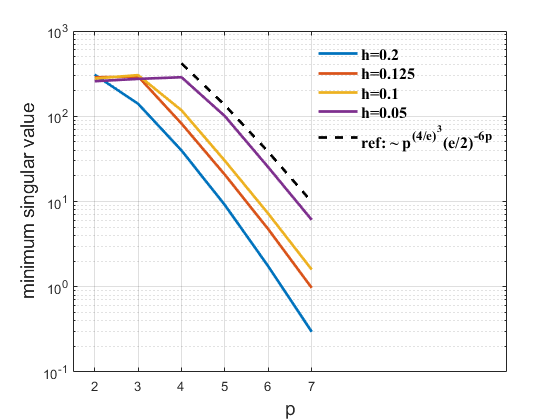}}
      \caption{The maximum/minimum singular value and the spectral condition number 
            of $\mathcal{K}(\mathbf{A}_{p-1})$ for a cube example, versus $h$ (at left) and 
            versus $p$ (at right).}
      \label{fig:3d_p_collocation}
\end{figure}

\begin{figure}[htbp]
      \centering  
      \subfigure[Condition number vs. $h$]{
            \label{subfig:3d_p_condA_h_add}
            \includegraphics[width=0.39\linewidth]{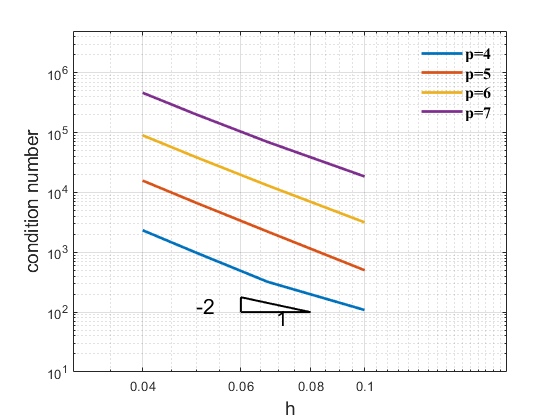}}    
      \subfigure[Condition number vs. $p$]{
            \label{subfig:3d_p_condA_p_add}
            \includegraphics[width=0.39\linewidth]{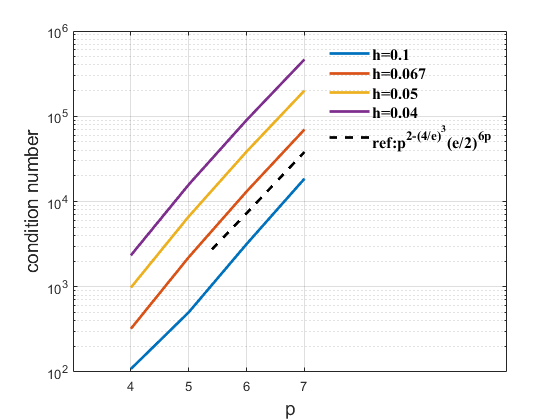}}\\             
      \subfigure[Maximum singular value vs. $h$]{
            \label{subfig:3d_p_maxA_h_add}
            \includegraphics[width=0.39\linewidth]{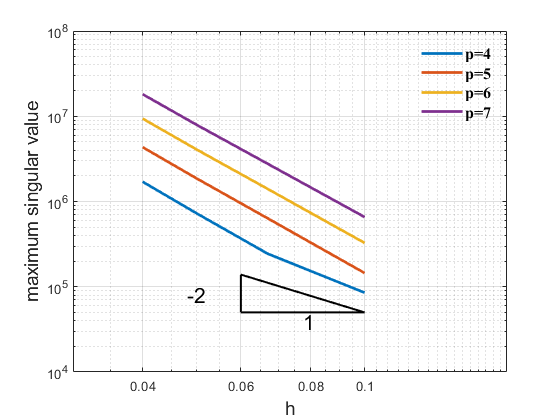}}                
      \subfigure[Maximum singular value vs. $p$]{
            \label{subfig:3d_p_maxA_p_add}
            \includegraphics[width=0.39\linewidth]{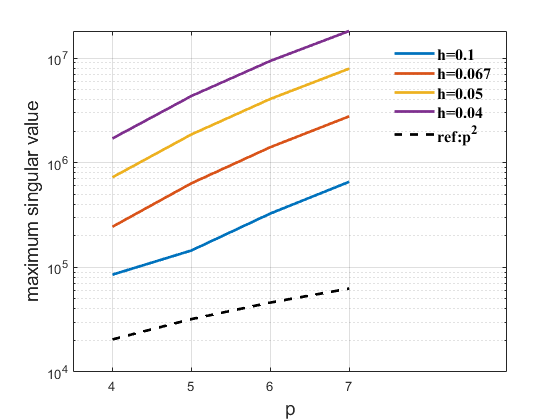}}\\
      \subfigure[Minimum singular value vs. $h$]{
            \label{subfig:3d_p_minA_h_add}
            \includegraphics[width=0.39\linewidth]{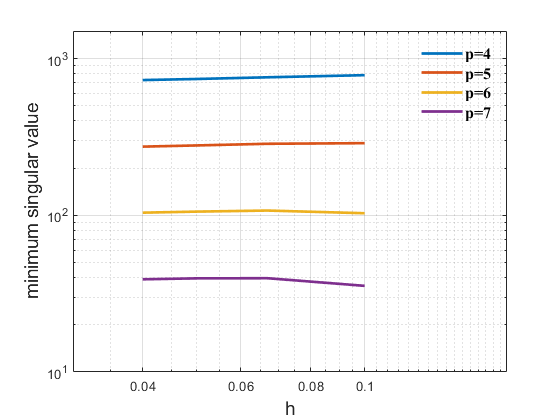}}               
      \subfigure[Minimum singular value vs. $p$]{
            \label{subfig:3d_p_minA_p_add}
            \includegraphics[width=0.39\linewidth]{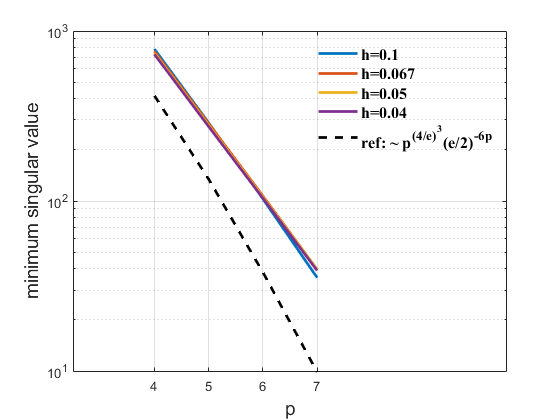}}
      \caption{The maximum/minimum singular value and the spectral condition number 
            of $\mathcal{K}(\mathbf{A}_{p-1})$ for $1/8$th of hollow sphere, 
            versus $h$ (at left) and versus $p$ (at right).}
      \label{fig:3d_p_collocation_add}
\end{figure}

In Figs.\ref{fig:1d_p_collocation}, \ref{fig:2d_p_collocation},
      \ref{fig:3d_p_collocation} and \ref{fig:3d_p_collocation_add},  
      we demonstrate the $\mathcal{K}(\mathbf{A}_{p-1})$,   
      $\sigma_{max}(\mathbf{A}_{p-1})$ 
      and $\sigma_{min}(\mathbf{A}_{p-1})$ versus 
      $h$ (at left) and $p$ (at right) for $d=1,2,3$.  
As shown in Figs. \ref{subfig:1d_p_condA_h}, \ref{subfig:2d_p_condA_h}, 
      \ref{subfig:3d_p_condA_h} and \ref{subfig:3d_p_condA_h_add}, 
      $\mathcal{K}(\mathbf{A}_{p-1})$ increases at a rate of about $h^{-2}$.
In Fig. \ref{subfig:3d_p_condA_h_add}, $\mathcal{K}(\mathbf{A}_{p-1})$ tends 
      to grow at a slightly higher rate than $h^{-2}$ in a case for small $h$.
Figs. \ref{subfig:1d_p_condA_p}, \ref{subfig:2d_p_condA_p}, 
      \ref{subfig:3d_p_condA_p} and \ref{subfig:3d_p_condA_p_add} 
      show that $\mathcal{K}(\mathbf{A}_{p-1})$ increases at a rate of 
      about $p^{2-(4/e)^{d}}(e/2)^{2dp}$, which is essentially in 
      line with our estimates.
Numerical results show that the condition number of $\mathbf{A}_{p-1}$ 
      with respect to $p$ basically satisfies 
      $\mathcal{K}(\mathbf{A}_{p-1})\thicksim h^{-2}p^{2}$ 
      when $h\lesssim e^{-pd/4}p^{d/8}$.
Additionally, when $h\gtrsim e^{-pd/4}p^{d/8}$, 
      the condition number grows exponentially with $p$ and 
      tends to be stable with $h-$refinement, 
      which is due to the increase in the minimum singular value at this time. 

As for the maximum singular value $\sigma_{max}(\mathbf{A}_{p-1})$,  
      Figs.\ref{subfig:1d_p_maxA_h}, \ref{subfig:2d_p_maxA_h},
      \ref{subfig:3d_p_maxA_h} and \ref{subfig:3d_p_maxA_h_add} 
      show that $\sigma_{max}(\mathbf{A}_{p-1})$ increases 
      at a rate of about $h^{-2}$ with $h-$refinement for any value of $p\geq 2$.
Figs.\ref{subfig:1d_p_maxA_p}, \ref{subfig:2d_p_maxA_p},
      \ref{subfig:3d_p_maxA_p} and \ref{subfig:3d_p_maxA_p_add} 
      illustrate a growth rate of about $p^{2}$ for any 
      value of $h>0$.
However, in Fig. \ref{subfig:3d_p_maxA_h_add}, 
      it can be observed that for the hollow sphere example, 
      the maximum singular value $\sigma_{max}(\mathbf{A}_{p-1})$ initially increases 
at a rate of about $h^{-2}$ with h-refinement. 
In a case for small $h$, $\sigma_{max}(\mathbf{A}_{p-1})$ exhibits a growth rate 
      slightly higher than $h^{-2}$. 
Moreover, Fig. \ref{subfig:3d_p_maxA_p_add} demonstrates that as $p$ increases, 
      the growth rate of $\sigma_{max}(\mathbf{A}_{p-1})$ is slightly 
      higher than $p^{2}$.

Regarding $\sigma_{min}(\mathbf{A}_{p-1})$, it initially increases by $h^{-2}$,   
      then stabilizes with further $h-$refinement.
As $p$ increases, the mesh size $h$ required for $\sigma_{min}(\mathbf{A}_{p-1})$ 
      to remain stable becomes smaller.
Fig.\ref{subfig:1d_p_minA_h} shows that for $p=12$, 
      the breakpoint for $h$ lies between $0.05$ and $0.1$, 
      while for $p=14$, it lies between $0.033$ and $0.05$.  
The numerical results about $\sigma_{min}(\mathbf{A}_{p-1})$ 
      show that $\sigma_{min}(\mathbf{A}_{p-1})$ remains independent of 
      $h$ and $p$ when $h\lesssim e^{-pd/4}p^{d/8}$. 
When $h$ is relatively large, 
      the support interval of the basis functions 
      can no longer be completely contained within the parameter domain.
This affects the linear independence of the column vectors of $\mathbf{A}_{p-1}$, 
      leading to an increase in $\sigma_{min}(\mathbf{A}_{p-1})$ with $h-$refinement 
      and a decrease of about $(e/2)^{-2dp}p^{(4/e)^{d}}$ as $p$ increases.
For the hollow sphere example, Fig. \ref{subfig:3d_p_minA_h_add} shows that 
      $\sigma_{min}(\mathbf{A}_{p-1})$ tends to be stable with $h$-refinement. 
Fig. \ref{subfig:3d_p_minA_p_add} shows as $p$ increases, 
      $\sigma_{min}(\mathbf{A}_{p-1})$ decreases by about 
      $(e/2)^{-6p}p^{(4/e)^{3}}$, which is basically consistent 
      with our estimates.

As shown in Fig. \ref{fig:3d_p_collocation_add},  
      the numerical results of the collocation matrix for the hollow sphere 
      example deviate slightly from our estimates, 
      which may be due to different geometric mappings.    
However, in general, the results are basically in line with our estimates, 
      which suggests that the spectral properties of the collocation matrix 
      behave approximately the same under different geometric mappings.

According to (\ref{eq:estimate_pA_cond}) and (\ref{eq:estimate_IGA-G_p}), 
      it can be observed that when $h\lesssim e^{-pd/2}$, 
      the stiffness matrix $\mathbf{K}_{p-1}$ in IGA-G is slightly better 
      conditioned than the collocation matrix in IGA-L.
However, this condition on $h$ is quite stringent in practice.
When $h\gtrsim e^{-pd/4}p^{d/8}$,
      the condition number of the collocation matrix satisfies 
      $\mathcal{K}(\mathbf{A}_{p-1})\thicksim (e/2)^{2dp}p^{2-(4/e)^{d}}$, 
      which is much better conditioned than $\mathbf{K}_{p-1}$.

\subsection{IGA-$C^{p-1}$ mass matrix} 
In this section, we study the spectral properties of the IGA-$C^{p-1}$ 
      mass matrix $\mathbf{M}_{p-1}$.
The mesh size $h$ and the degree $p$ are set to be
      the same as those used in the previous section.
The numerical results show that for $h>0$, $p\geq 2$ and $d=1,2,3$, 
      $\sigma_{min}(\mathbf{M}_{p-1})$, $\sigma_{max}(\mathbf{M}_{p-1})$ and 
      $\mathcal{K}(\mathbf{M}_{p-1})$ behave as:
\begin{align}
      \label{eq:estimate_pM_min}
            &\sigma_{min}(\mathbf{M}_{p-1})\thicksim e^{-dp/2}\\
            \label{eq:estimate_pM_max}
            &\sigma_{max}(\mathbf{M}_{p-1})\thicksim c_{2}, \\
            \label{eq:estimate_pM_cond}
            &\mathcal{K}(\mathbf{M}_{p-1})\thicksim e^{dp/2} 
\end{align}
where $c_{2}$ is independent of $h$ and $p$.
      
      \begin{figure}[htbp]
            \centering  
            \subfigure[Condition number vs. $h$]{
                  \label{subfig:1d_p_condM_h}
                  \includegraphics[width=0.39\linewidth]{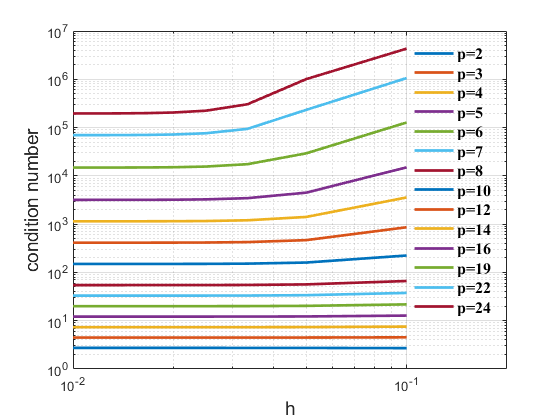}}    
            \subfigure[Condition number vs. $p$]{
                  \label{subfig:1d_p_condM_p}
                  \includegraphics[width=0.39\linewidth]{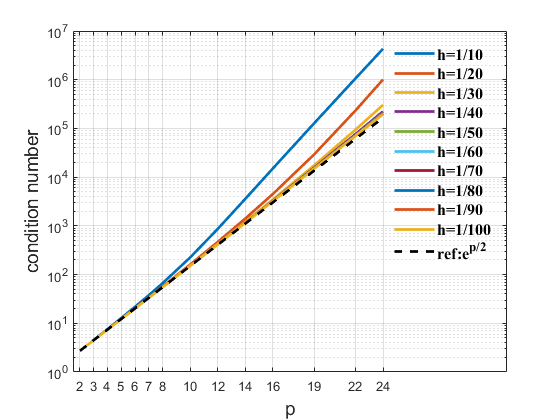}}\\             
            \subfigure[Maximum singular value vs. $h$]{
                  \label{subfig:1d_p_maxM_h}
                  \includegraphics[width=0.39\linewidth]{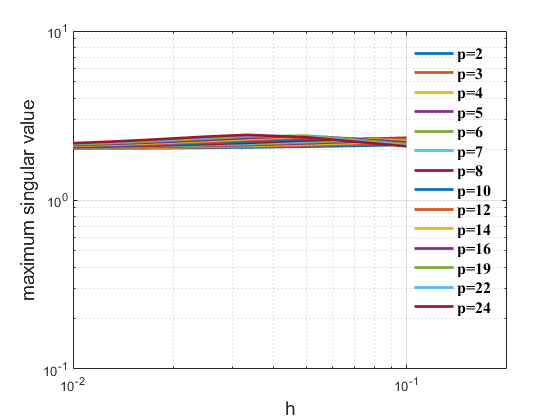}}                
            \subfigure[Maximum singular value vs. $p$]{
                  \label{subfig:1d_p_maxM_p}
                  \includegraphics[width=0.39\linewidth]{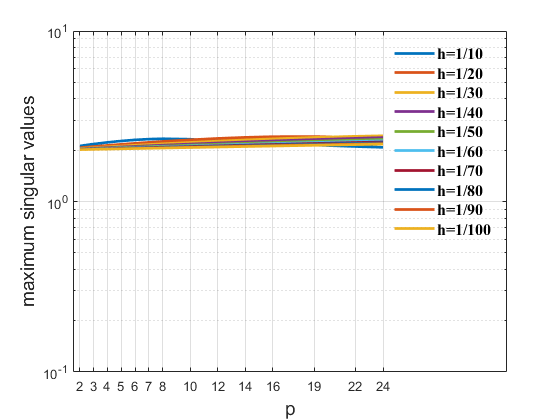}}\\
            \subfigure[Minimum singular value vs. $h$]{
                  \label{subfig:1d_p_minM_h}
                  \includegraphics[width=0.39\linewidth]{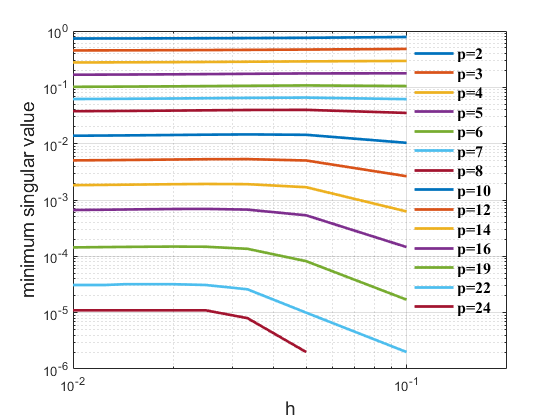}}               
            \subfigure[Minimum singular value vs. $p$]{
                  \label{subfig:1d_p_minM_p}
                  \includegraphics[width=0.39\linewidth]{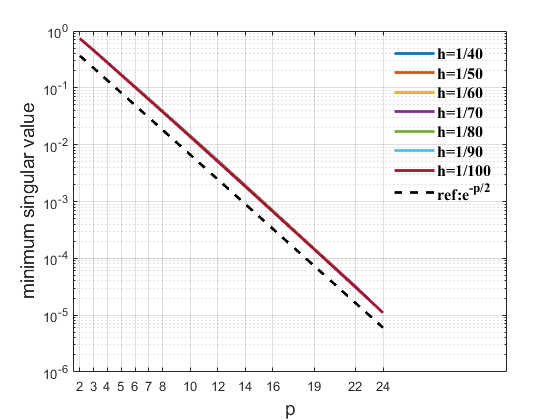}}
            \caption{The maximum/minimum singular value and the spectral condition number 
            of $\mathcal{K}(\mathbf{M}_{p-1})$ for $d=1$, versus $h$ (at left) and 
            versus $p$ (at right).}
            \label{fig:1d_p_mass}
      \end{figure}

      \begin{figure}[htbp]
            \centering  
            \subfigure[Condition number vs. $h$]{
                  \label{subfig:2d_p_condM_h}
                  \includegraphics[width=0.39\linewidth]{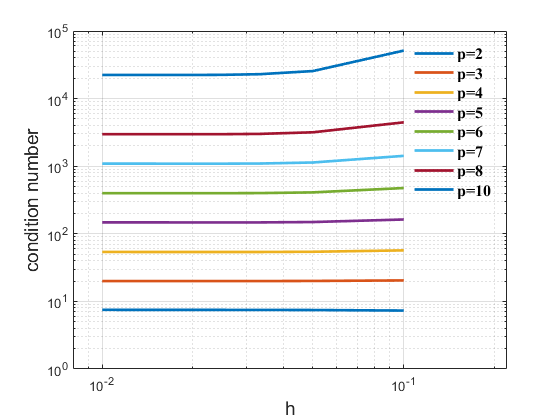}}    
            \subfigure[Condition number vs. $p$]{
                  \label{subfig:2d_p_condM_p}
                  \includegraphics[width=0.39\linewidth]{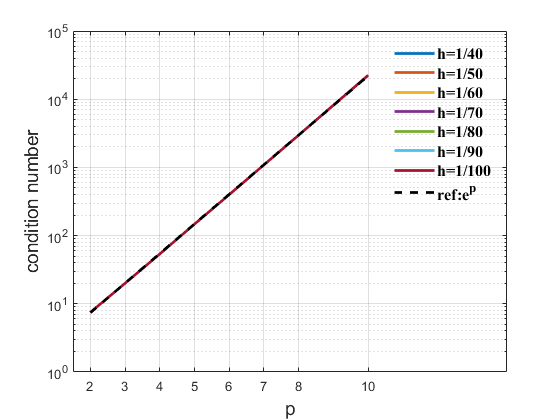}}\\             
            \subfigure[Maximum singular value vs. $h$]{
                  \label{subfig:2d_p_maxM_h}
                  \includegraphics[width=0.39\linewidth]{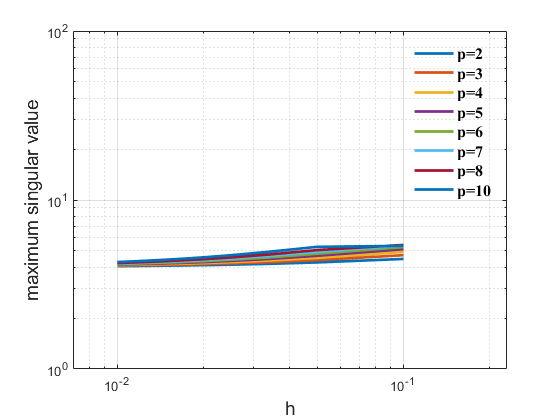}}                
            \subfigure[Maximum singular value vs. $p$]{
                  \label{subfig:2d_p_maxM_p}
                  \includegraphics[width=0.39\linewidth]{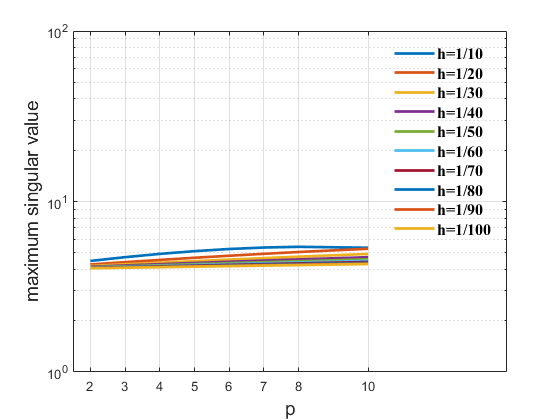}}\\
            \subfigure[Minimum singular value vs. $h$]{
                  \label{subfig:2d_p_minM_h}
                  \includegraphics[width=0.39\linewidth]{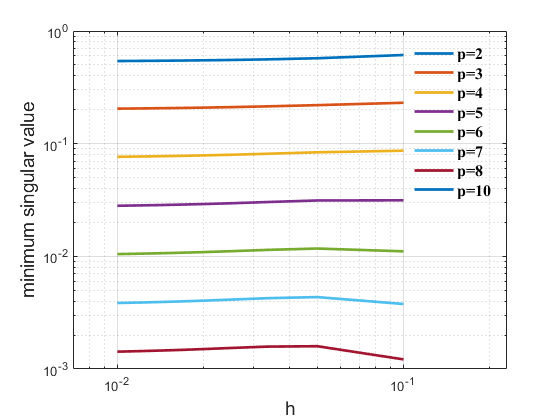}}               
            \subfigure[Minimum singular value vs. $p$]{
                  \label{subfig:2d_p_minM_p}
                  \includegraphics[width=0.39\linewidth]{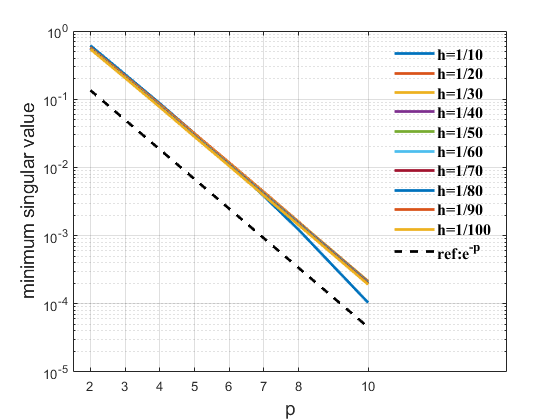}}
            \caption{The maximum/minimum singular value and the spectral condition number 
            of $\mathcal{K}(\mathbf{M}_{p-1})$ for $d=2$, versus $h$ (at left) and 
            versus $p$ (at right).}
            \label{fig:2d_p_mass}
      \end{figure}

      \begin{figure}[htbp]
            \centering  
            \subfigure[Condition number vs. $h$]{
                  \label{subfig:3d_p_condM_h}
                  \includegraphics[width=0.38\linewidth]{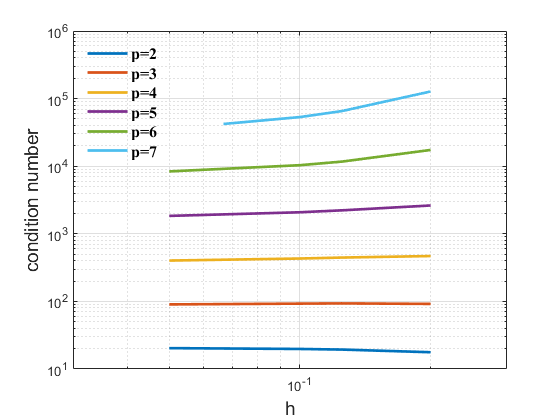}}    
            \subfigure[Condition number vs. $p$]{
                  \label{subfig:3d_p_condM_p}
                  \includegraphics[width=0.38\linewidth]{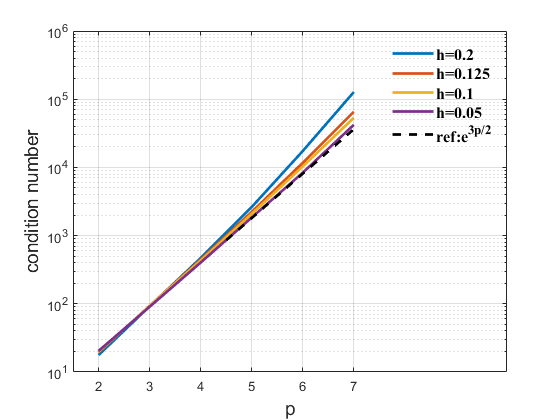}}\\             
            \subfigure[Maximum singular value vs. $h$]{
                  \label{subfig:3d_p_maxM_h}
                  \includegraphics[width=0.38\linewidth]{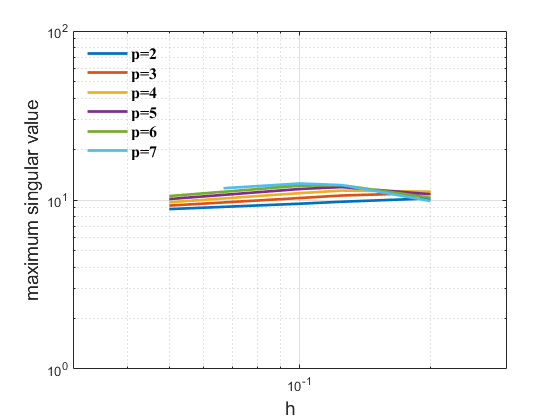}}                
            \subfigure[Maximum singular value vs. $p$]{
                  \label{subfig:3d_p_maxM_p}
                  \includegraphics[width=0.38\linewidth]{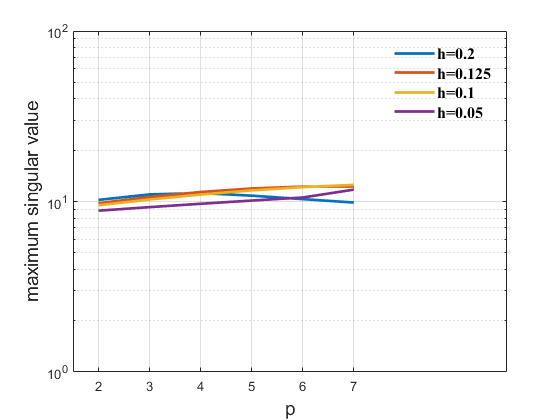}}\\
            \subfigure[Minimum singular value vs. $h$]{
                  \label{subfig:3d_p_minM_h}
                  \includegraphics[width=0.38\linewidth]{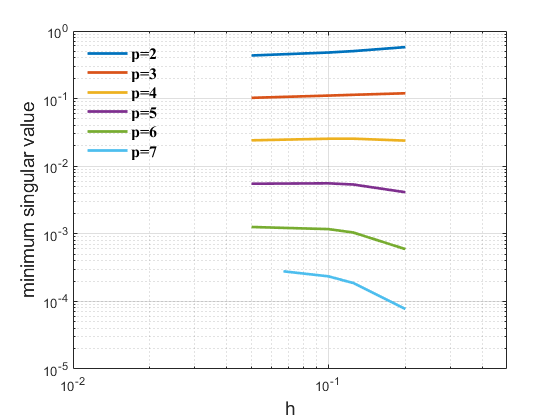}}               
            \subfigure[Minimum singular value vs. $p$]{
                  \label{subfig:3d_p_minM_p}
                  \includegraphics[width=0.38\linewidth]{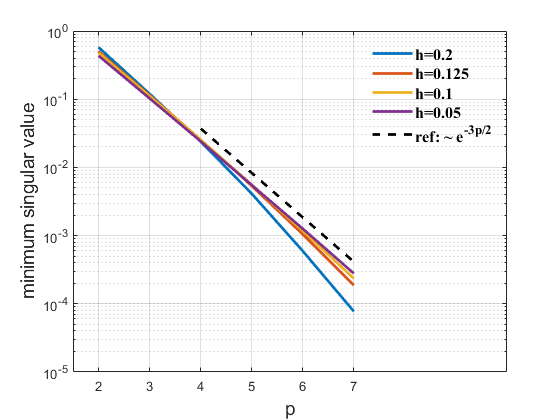}}
            \caption{The maximum/minimum singular value and the spectral condition number 
                  of $\mathcal{K}(\mathbf{M}_{p-1})$ for a cube example, versus $h$ (at left) and 
                  versus $p$ (at right).}
            \label{fig:3d_p_mass}
      \end{figure} 

      \begin{figure}[htbp]
            \centering  
            \subfigure[Condition number vs. $h$]{
                  \label{subfig:3d_p_condM_h_add}
                  \includegraphics[width=0.38\linewidth]{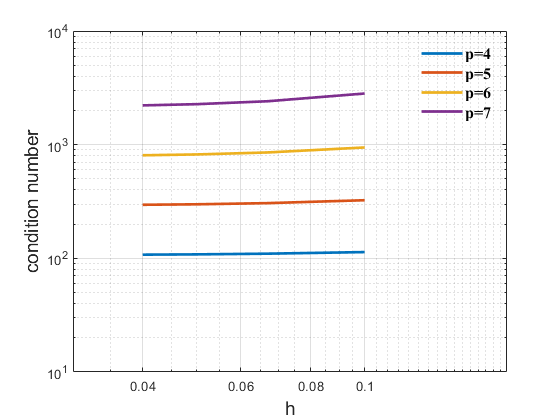}}    
            \subfigure[Condition number vs. $p$]{
                  \label{subfig:3d_p_condM_p_add}
                  \includegraphics[width=0.38\linewidth]{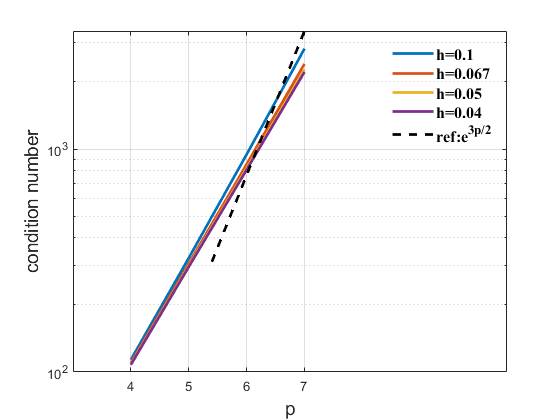}}\\             
            \subfigure[Maximum singular value vs. $h$]{
                  \label{subfig:3d_p_maxM_h_add}
                  \includegraphics[width=0.38\linewidth]{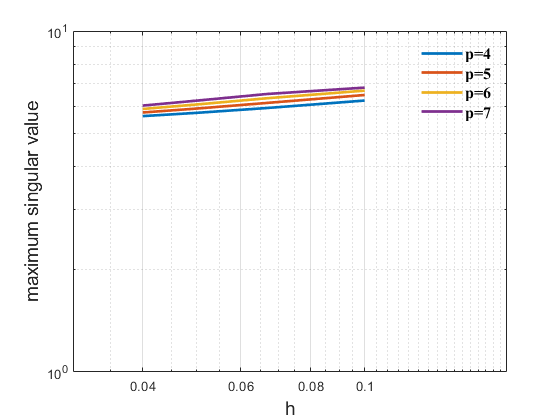}}                
            \subfigure[Maximum singular value vs. $p$]{
                  \label{subfig:3d_p_maxM_p_add}
                  \includegraphics[width=0.38\linewidth]{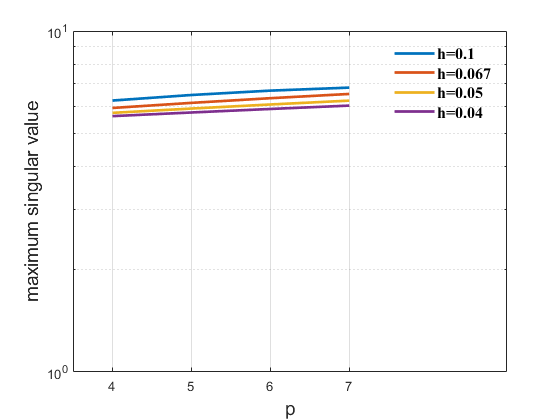}}\\
            \subfigure[Minimum singular value vs. $h$]{
                  \label{subfig:3d_p_minM_h_add}
                  \includegraphics[width=0.38\linewidth]{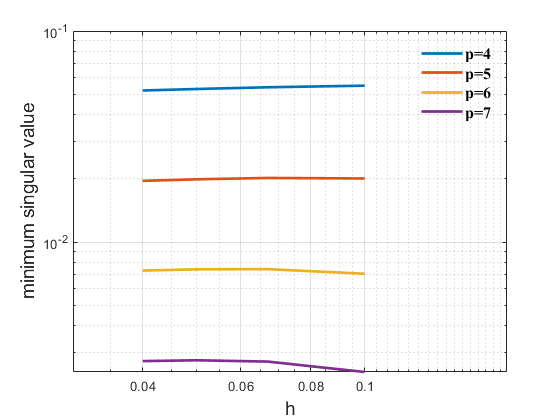}}               
            \subfigure[Minimum singular value vs. $p$]{
                  \label{subfig:3d_p_minM_p_add}
                  \includegraphics[width=0.38\linewidth]{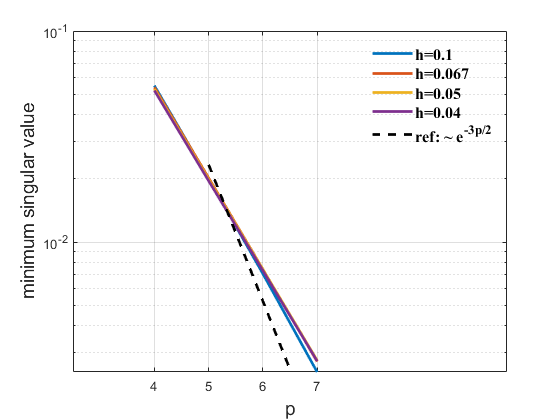}}
            \caption{The maximum/minimum singular value and the spectral condition number 
                  of $\mathcal{K}(\mathbf{M}_{p-1})$ for $1/8$th of hollow sphere, 
                  versus $h$ (at left) and versus $p$ (at right).}
            \label{fig:3d_p_mass_add}
      \end{figure}

In Figs.\ref{fig:1d_p_mass}, \ref{fig:2d_p_mass}, 
      \ref{fig:3d_p_mass} and \ref{fig:3d_p_mass_add}, 
      we show the $\mathcal{K}(\mathbf{M}_{p-1})$,   
      $\sigma_{max}(\mathbf{M}_{p-1})$ and $\sigma_{min}(\mathbf{M}_{p-1})$ versus 
      $h$ (at left) and $p$ (at right) for $d=1,2,3$.    
For the condition number $\mathcal{K}(\mathbf{M}_{p-1})$,
Figs.\ref{subfig:1d_p_condM_h}, 
      \ref{subfig:2d_p_condM_h}, \ref{subfig:3d_p_condM_h} and 
      \ref{subfig:3d_p_condM_h_add}   
demonstrate that $\mathcal{K}(\mathbf{M}_{p-1})$ does not vary with $h$ 
      when $h$ is relatively small.
Figs.\ref{subfig:1d_p_condM_p}, \ref{subfig:2d_p_condM_p} 
      and \ref{subfig:3d_p_condM_p} 
      show that for $d=1$, 
      $\mathcal{K}(\mathbf{M}_{p-1})$ increases by about $e^{p/2}$; 
      for $d=2$, it increases as $e^{p}$;
      and for $d=3$, it increases by about $e^{3p/2}$.
However, Fig. \ref{subfig:3d_p_condM_p_add} demonstrates that 
      for the hollow sphere example, 
      the growth rate of $\mathcal{K}(\mathbf{M}_{p-1})$ is 
      slightly lower than $e^{3p/2}$.
For the maximum singular value $\sigma_{max}(\mathbf{M}_{p-1})$, 
      numerical results
      demonstrate that $\sigma_{max}(\mathbf{M}_{p-1})$ 
      is independent of both $h$ and $p$ for $h>0$ and $p\geq 2$.
This is due to the relatively uniform distribution of the 
      columns of the mass matrix formed by 
      the internal NURBS basis functions. 

As shown in Figs.\ref{subfig:1d_p_minM_h}, \ref{subfig:2d_p_minM_h}, 
      \ref{subfig:3d_p_minM_h} and \ref{subfig:3d_p_minM_h_add}, 
      the minimum singular value of the mass matrix $\sigma_{min}(\mathbf{M}_{p-1})$ 
      essentially does not change with $h-$refinement.
In Fig.\ref{subfig:1d_p_minM_h}, it can be observed that 
      for relatively large degrees $p$, 
      $\sigma_{min}(\mathbf{M}_{p-1})$ initially 
      increases with $h-$refinement, and then remains constant.
Figs.\ref{subfig:1d_p_minM_p}, \ref{subfig:2d_p_minM_p}, 
      \ref{subfig:3d_p_minM_p} and \ref{subfig:3d_p_minM_p_add} show that 
      as $p$ increases, 
      $\sigma_{min}(\mathbf{M}_{p-1})$ decreases 
      exponentially with $p$.
When $d=1$, 
      $\sigma_{min}(\mathbf{M}_{p-1})$ decreases by about $e^{-p/2}$ as $p$; 
      when $d=2$, it decreases as $e^{-p}$;
      and when $d=3$, it decreases by about $e^{-3p/2}$.
      Moreover, Fig. \ref{subfig:3d_p_minM_p_add} shows that 
      for the hollow sphere example, 
      $\sigma_{min}(\mathbf{M}_{p-1})$ decreases slightly 
      slower than $e^{-3p/2}$.
In general, the numerical results shown in Figs.\ref{fig:1d_p_mass}, 
      \ref{fig:2d_p_mass}, \ref{fig:3d_p_mass} and \ref{fig:3d_p_mass_add}  
      are basically in line with our estimates (\ref{eq:estimate_pM_min})-
(\ref{eq:estimate_pM_cond}).
Compared with the condition number estimations of $\mathbf{M}_{p-1}$ 
      in IGA-G (\ref{eq:estimate_IGA-G_p}), 
      we find that $\mathbf{M}_{p-1}$ in IGA-L is better conditioned for any value $h>0$. 
\subsection{IGA-$C^{1}$ collocation matrix} 
In this section, we study the spectral properties of the IGA-$C^{1}$ 
      collocation matrix $\mathbf{A}_{1}$.
The mesh size $h$ varies from $0.1$ to $0.01$ for $d=1$.
For $d=2$, $h$ is set to decrease from $0.2$ to $0.05$.
For $d=3$, $h$ is set to decrease from $0.2$ to $0.1$.
The degree $p$ starts at $2$, 
      and increases to $16,7,4$ for $d=1,2,3$, respectively. 

The numerical results show that for $h>0$, $p\geq 2$ and $d=1,2,3$, 
      $\sigma_{min}(\mathbf{A}_{1})$, $\sigma_{max}(\mathbf{A}_{1})$ and 
      $\mathcal{K}(\mathbf{A}_{1})$ behave as:
\begin{align}
      \label{eq:estimate_1A_min}
            &\sigma_{min}(\mathbf{A}_{1})\thicksim c_{3}\\
            \label{eq:estimate_1A_max}
            &\sigma_{max}(\mathbf{A}_{1})\thicksim h^{-2}p^{5/2}\\
            \label{eq:estimate_1A_cond}
            &\mathcal{K}(\mathbf{A}_{1})\thicksim h^{-2}p^{5/2} 
\end{align}
where $c_{3}$ is independent of $h$ and $p$.

\begin{figure}[htbp]   
      \centering  
      \subfigure[Condition number vs. $h$]{
            \label{subfig:1d_1_condA_h}
            \includegraphics[width=0.39\linewidth]{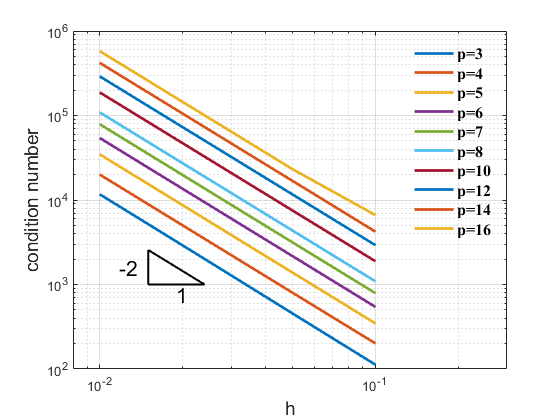}}    
      \subfigure[Condition number vs. $p$]{
            \label{subfig:1d_1_condA_p}
            \includegraphics[width=0.39\linewidth]{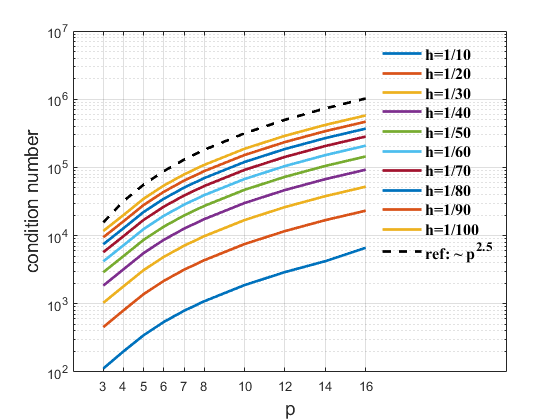}}\\             
      \subfigure[Maximum singular value vs. $h$]{
            \label{subfig:1d_1_maxA_h}
            \includegraphics[width=0.39\linewidth]{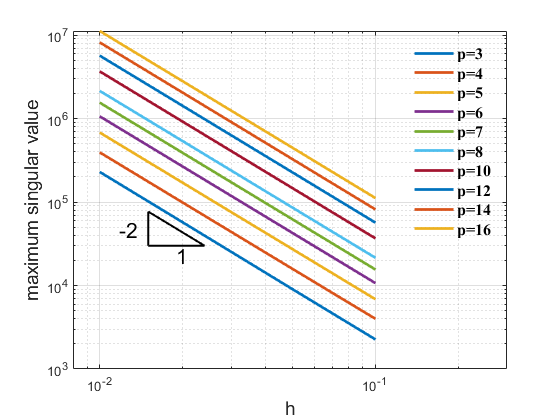}}                
      \subfigure[Maximum singular value vs. $p$]{
            \label{subfig:1d_1_maxA_p}
            \includegraphics[width=0.39\linewidth]{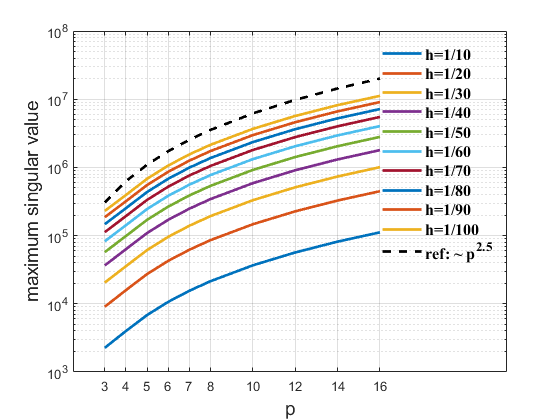}}\\
      \subfigure[Minimum singular value vs. $h$]{ 
            \label{subfig:1d_1_minA_h}
            \includegraphics[width=0.39\linewidth]{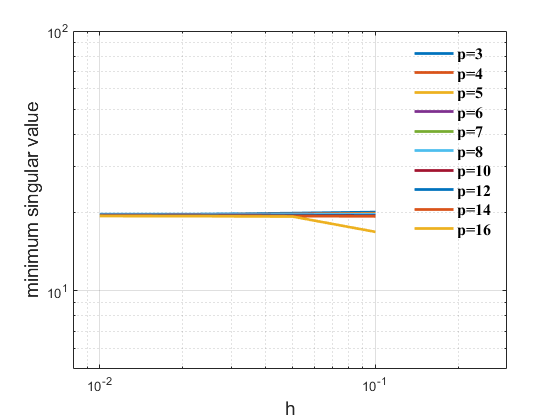}}               
      \subfigure[Minimum singular value vs. $p$]{
            \label{subfig:1d_1_minA_p}
            \includegraphics[width=0.39\linewidth]{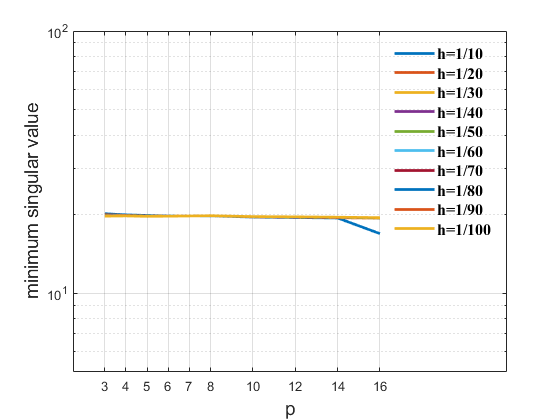}}
      \caption{The maximum/minimum singular value and the spectral condition number 
      of $\mathcal{K}(\mathbf{A}_{1})$ for $d=1$, versus $h$ (at left) and 
      versus $p$ (at right).}  
      \label{fig:1d_1_collocation}   
\end{figure}

\begin{figure}[htbp]
      \centering  
      \subfigure[Condition number vs. $h$]{
            \label{subfig:2d_1_condA_h}
            \includegraphics[width=0.39\linewidth]{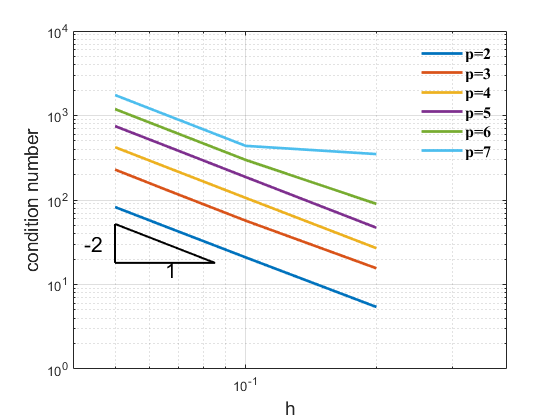}}    
      \subfigure[Condition number vs. $p$]{
            \label{subfig:2d_1_condA_p}
            \includegraphics[width=0.39\linewidth]{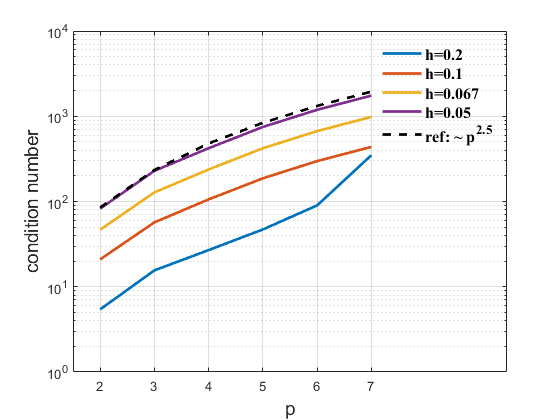}}\\             
      \subfigure[Maximum singular value vs. $h$]{
            \label{subfig:2d_1_maxA_h}
            \includegraphics[width=0.39\linewidth]{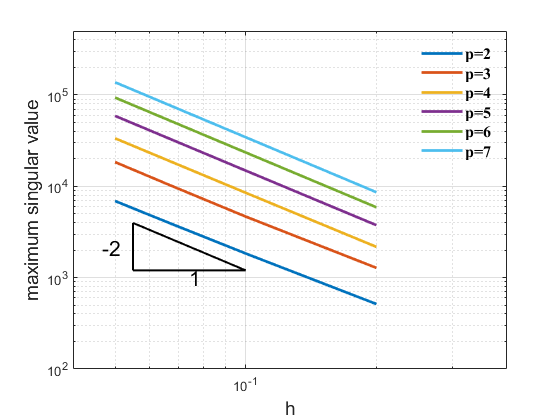}}                
      \subfigure[Maximum singular value vs. $p$]{
            \label{subfig:2d_1_maxA_p}
            \includegraphics[width=0.39\linewidth]{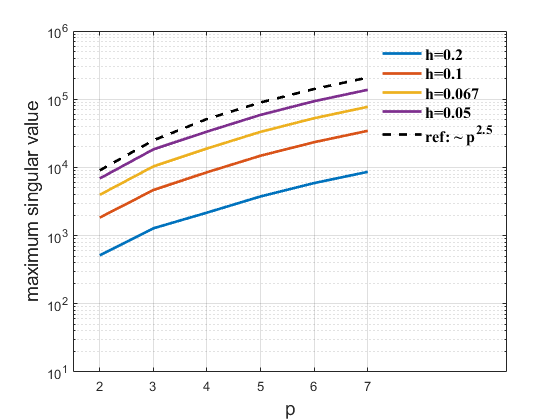}}\\
      \subfigure[Minimum singular value vs. $h$]{
            \label{subfig:2d_1_minA_h}
            \includegraphics[width=0.39\linewidth]{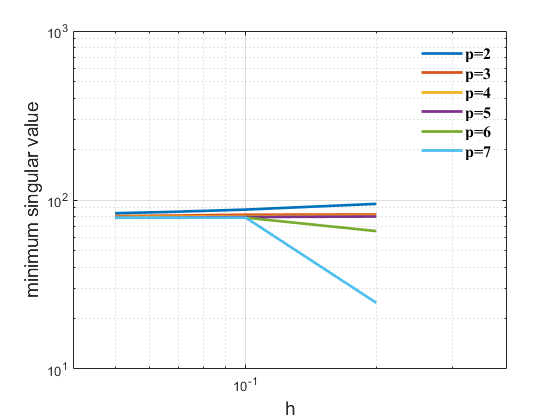}}               
      \subfigure[Minimum singular value vs. $p$]{
            \label{subfig:2d_1_minA_p}
            \includegraphics[width=0.39\linewidth]{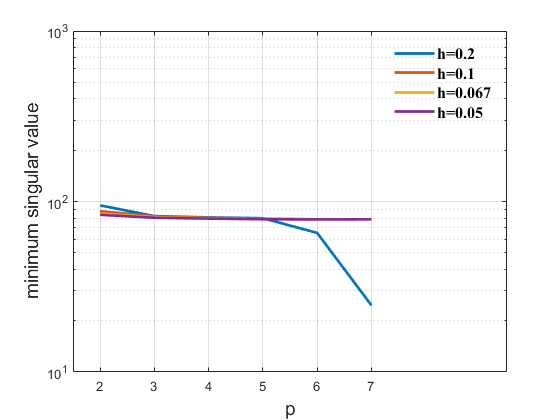}}
      \caption{The maximum/minimum singular value and the spectral condition number 
      of $\mathcal{K}(\mathbf{A}_{1})$ for $d=2$, versus $h$ (at left) and 
      versus $p$ (at right).}
      \label{fig:2d_1_collocation}
\end{figure}

\begin{figure}[htbp]
      \centering  
      \subfigure[Condition number vs. $h$]{
            \label{subfig:3d_1_condA_h}
            \includegraphics[width=0.39\linewidth]{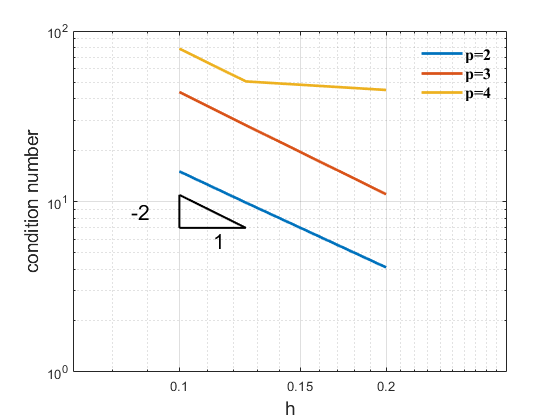}}    
      \subfigure[Condition number vs. $p$]{
            \label{subfig:3d_1_condA_p}
            \includegraphics[width=0.39\linewidth]{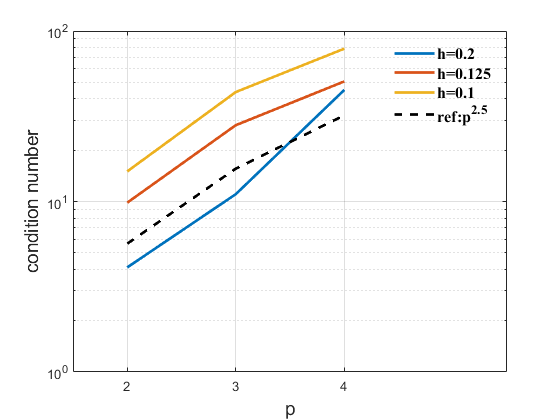}}\\             
      \subfigure[Maximum singular value vs. $h$]{
            \label{subfig:3d_1_maxA_h}
            \includegraphics[width=0.39\linewidth]{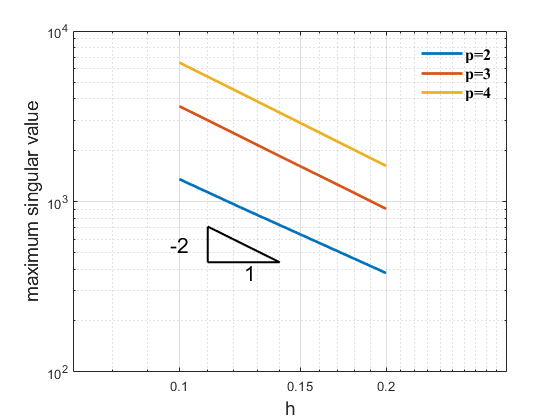}}                
      \subfigure[Maximum singular value vs. $p$]{
            \label{subfig:3d_1_maxA_p}
            \includegraphics[width=0.39\linewidth]{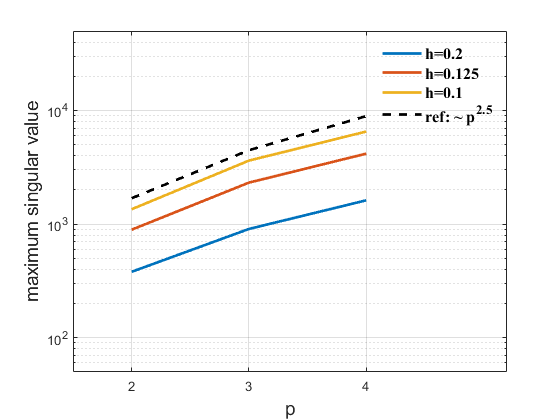}}\\
      \subfigure[Minimum singular value vs. $h$]{
            \label{subfig:3d_1_minA_h}
            \includegraphics[width=0.39\linewidth]{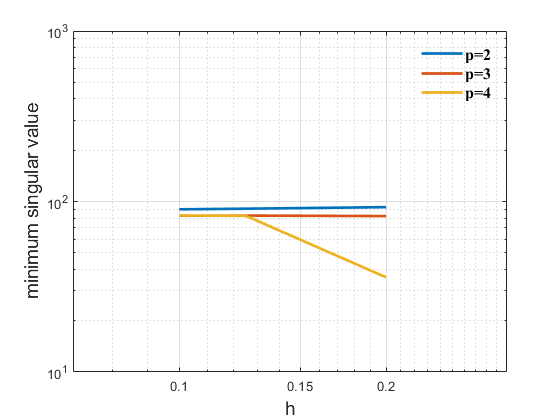}}               
      \subfigure[Minimum singular value vs. $p$]{
            \label{subfig:3d_1_minA_p}
            \includegraphics[width=0.39\linewidth]{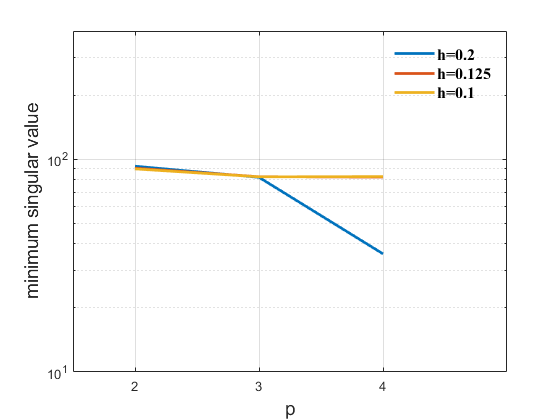}}
      \caption{The maximum/minimum singular value and the spectral condition number 
      of $\mathcal{K}(\mathbf{A}_{1})$ for $d=3$, versus $h$ (at left) and 
      versus $p$ (at right).}
      \label{fig:3d_1_collocation}
\end{figure}

Figs.\ref{subfig:1d_1_maxA_h}, \ref{subfig:2d_1_maxA_h} and 
      \ref{subfig:3d_1_maxA_h} demonstrate that $\sigma_{max}(\mathbf{A}_{1})$ 
      grows as $h^{-2}$ with $h-$refinement. 
Moreover, Figs.\ref{subfig:1d_1_maxA_p}, \ref{subfig:2d_1_maxA_p} and 
      \ref{subfig:3d_1_maxA_p} show that $\sigma_{max}(\mathbf{A}_{1})$ grows 
      approximately as $p^{2.5}$ with $p-$refinement. 
Thus, we obtain the estimate $\sigma_{max}(\mathbf{A}_{1})\thicksim h^{-2}p^{2.5}$.

Figs.\ref{subfig:1d_1_minA_h}, \ref{subfig:2d_1_minA_h} and 
      \ref{subfig:3d_1_minA_h} show that $\sigma_{min}(\mathbf{A}_{1})$ 
      is independent of $h$ for any value $p\geq 2$.
Additionally, Figs.\ref{subfig:1d_1_minA_p}, \ref{subfig:2d_1_minA_p} and 
      \ref{subfig:3d_1_minA_p} demonstrate that $\sigma_{min}(\mathbf{A}_{1})$ 
      is also independent of $p$.
Thus, we obtain the estimate of the condition number 
      $\mathcal{K}(\mathbf{A}_{1})\thicksim h^{-2}p^{2.5}$.
The numerical results shown in Figs.\ref{fig:1d_1_collocation}, 
      \ref{fig:2d_1_collocation} and \ref{fig:3d_1_collocation} 
      confirm our estimations (\ref{eq:estimate_1A_min})-
      (\ref{eq:estimate_1A_cond}).

\subsection{IGA-$C^{1}$ mass matrix}
In this section, we study the spectral properties of the IGA-$C^{1}$ 
      mass matrix $\mathbf{M}_{1}$.
The mesh size $h$ and the degree $p$ are set to be
      the same as those used in the previous section.
The numerical results show that for $h>0$, $p\geq 2$ and $d=1,2,3$, 
      $\sigma_{min}(\mathbf{M}_{1})$, $\sigma_{max}(\mathbf{M}_{1})$ and 
      $\mathcal{K}(\mathbf{M}_{1})$ behave as:
\begin{align}
      \label{eq:estimate_1M_min}
            &\sigma_{min}(\mathbf{M}_{1})\thicksim (2d)^{-p}\\
            \label{eq:estimate_1M_max}
            &\sigma_{max}(\mathbf{M}_{1})\thicksim c_{4}\\
            \label{eq:estimate_1M_cond}
            &\mathcal{K}(\mathbf{M}_{1})\thicksim (2d)^{p}
\end{align}
where $c_{4}$ is independent of $h$ and $p$.

\begin{figure}[htbp]
      \centering  
      \subfigure[Condition number vs. $h$]{
            \label{subfig:1d_1_condM_h}
            \includegraphics[width=0.39\linewidth]{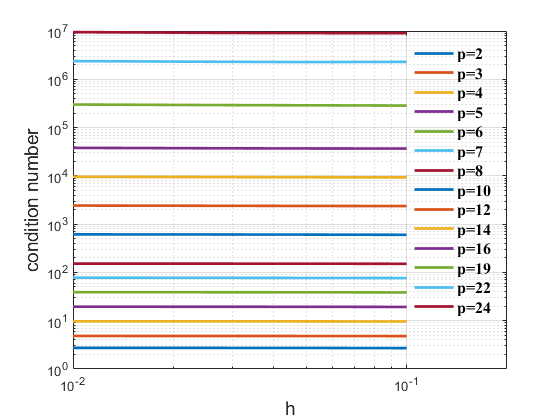}}    
      \subfigure[Condition number vs. $p$]{
            \label{subfig:1d_1_condM_p}
            \includegraphics[width=0.39\linewidth]{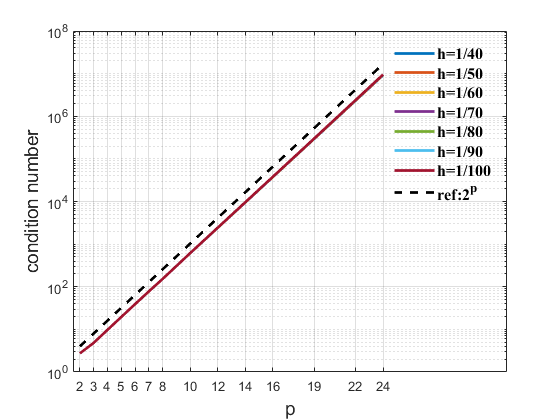}}\\             
      \subfigure[Maximum singular value vs. $h$]{
            \label{subfig:1d_1_maxM_h}
            \includegraphics[width=0.39\linewidth]{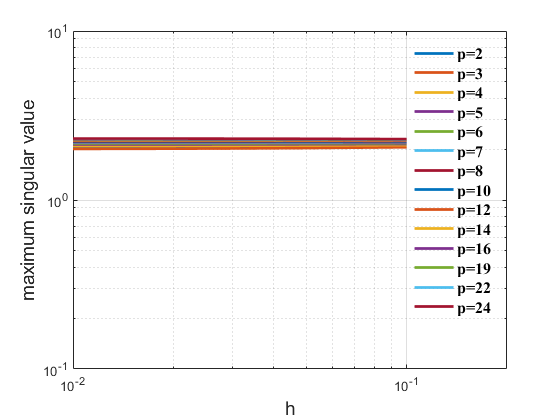}}                
      \subfigure[Maximum singular value vs. $p$]{
            \label{subfig:1d_1_maxM_p}
            \includegraphics[width=0.39\linewidth]{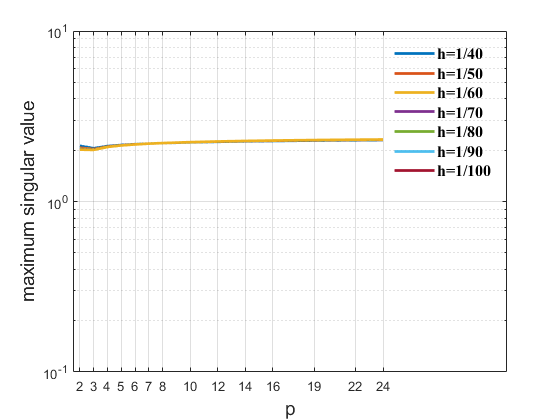}}\\
      \subfigure[Minimum singular value vs. $h$]{
            \label{subfig:1d_1_minM_h}
            \includegraphics[width=0.39\linewidth]{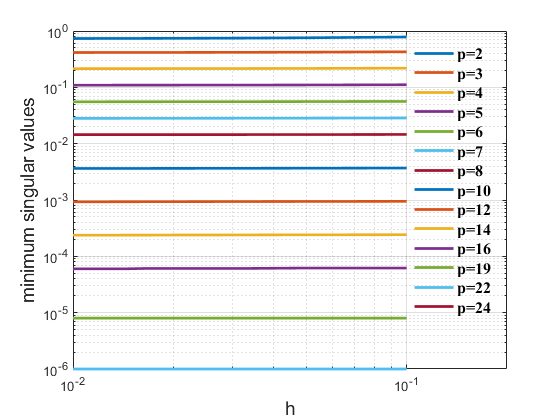}}               
      \subfigure[Minimum singular value vs. $p$]{
            \label{subfig:1d_1_minM_p}
            \includegraphics[width=0.39\linewidth]{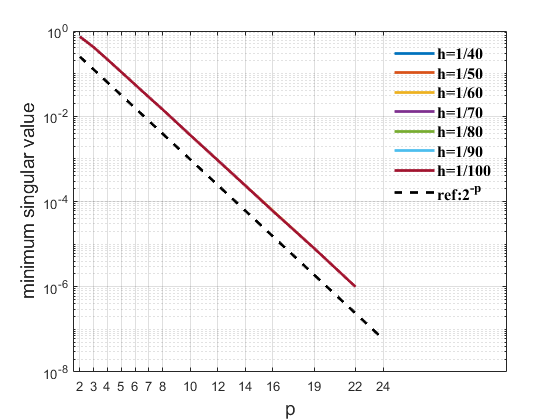}}
      \caption{The maximum/minimum singular value and the spectral condition number 
      of $\mathcal{K}(\mathbf{M}_{1})$ for $d=1$, versus $h$ (at left) and 
      versus $p$ (at right).}
      \label{fig:1d_1_mass}
\end{figure}

\begin{figure}[htbp]
      \centering  
      \subfigure[Condition number vs. $h$]{
            \label{subfig:2d_1_condM_h}
            \includegraphics[width=0.39\linewidth]{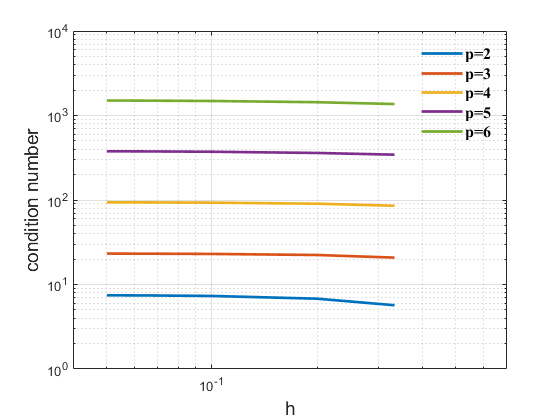}}    
      \subfigure[Condition number vs. $p$]{
            \label{subfig:2d_1_condM_p}
            \includegraphics[width=0.39\linewidth]{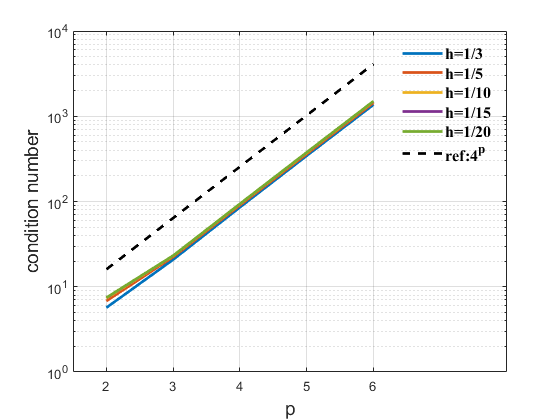}}\\             
      \subfigure[Maximum singular value vs. $h$]{
            \label{subfig:2d_1_maxM_h}
            \includegraphics[width=0.39\linewidth]{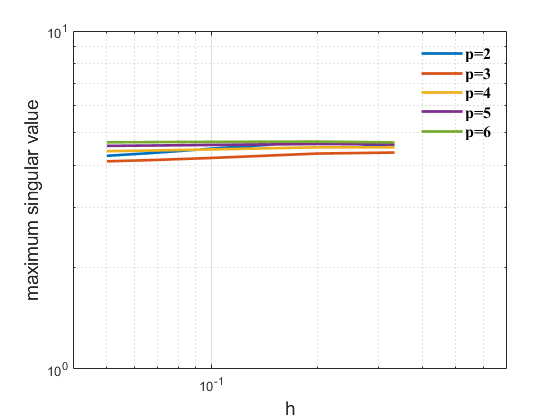}}                
      \subfigure[Maximum singular value vs. $p$]{
            \label{subfig:2d_1_maxM_p}
            \includegraphics[width=0.39\linewidth]{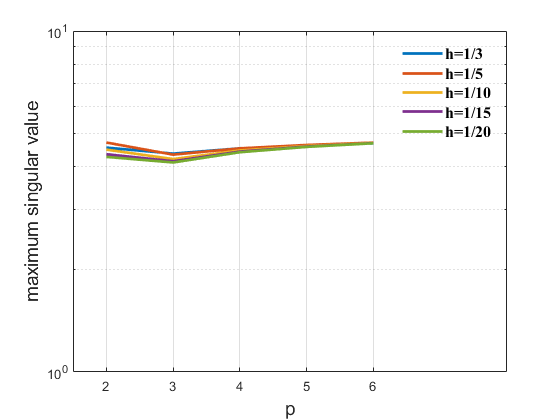}}\\
      \subfigure[Minimum singular value vs. $h$]{
            \label{subfig:2d_1_minM_h}
            \includegraphics[width=0.39\linewidth]{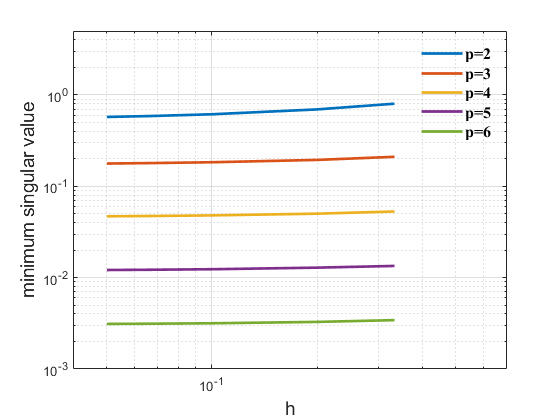}}               
      \subfigure[Minimum singular value vs. $p$]{
            \label{subfig:2d_1_minM_p}
            \includegraphics[width=0.39\linewidth]{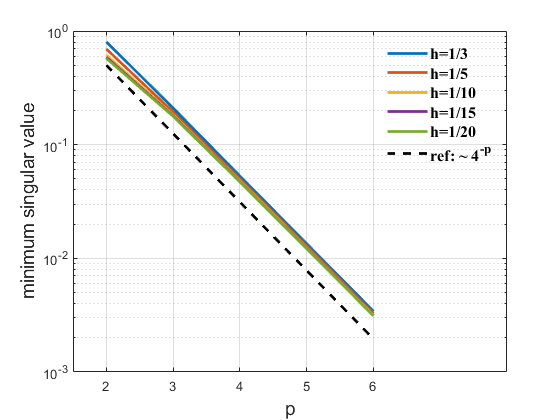}}
      \caption{The maximum/minimum singular value and the spectral condition number 
      of $\mathcal{K}(\mathbf{M}_{1})$ for $d=2$, versus $h$ (at left) and 
      versus $p$ (at right).}
      \label{fig:2d_1_mass}
\end{figure}

\begin{figure}[htbp]
      \centering  
      \subfigure[Condition number vs. $h$]{
            \label{subfig:3d_1_condM_h}
            \includegraphics[width=0.39\linewidth]{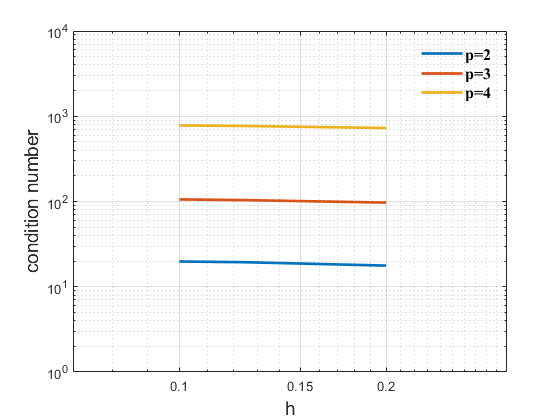}}    
      \subfigure[Condition number vs. $p$]{
            \label{subfig:3d_1_condM_p}
            \includegraphics[width=0.39\linewidth]{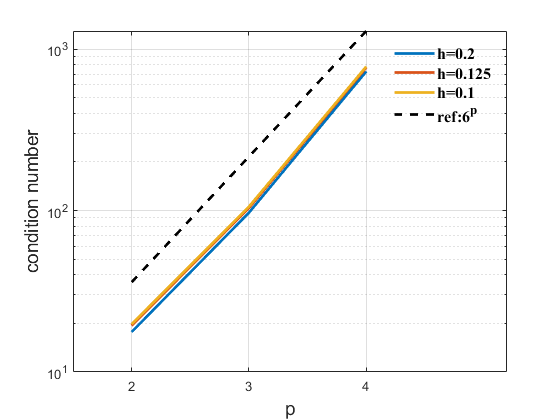}}\\             
      \subfigure[Maximum singular value vs. $h$]{
            \label{subfig:3d_1_maxM_h}
            \includegraphics[width=0.39\linewidth]{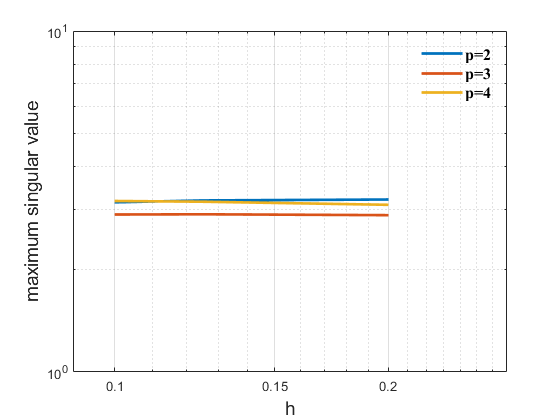}}                
      \subfigure[Maximum singular value vs. $p$]{
            \label{subfig:3d_1_maxM_p}
            \includegraphics[width=0.39\linewidth]{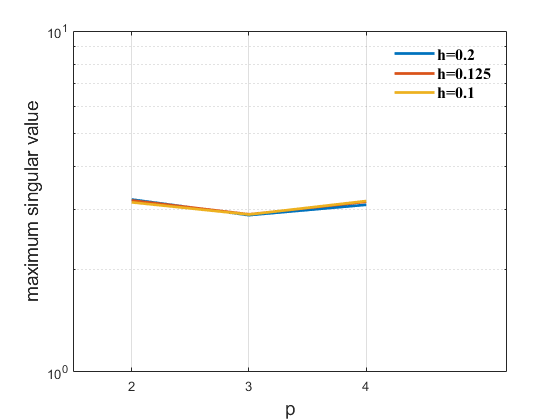}}\\
      \subfigure[Minimum singular value vs. $h$]{
            \label{subfig:3d_1_minM_h}
            \includegraphics[width=0.39\linewidth]{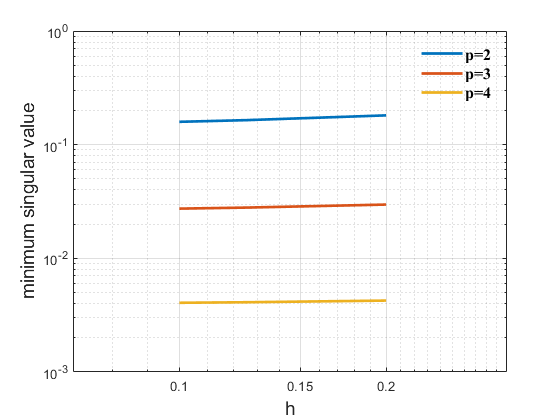}}               
      \subfigure[Minimum singular value vs. $p$]{
            \label{subfig:3d_1_minM_p}
            \includegraphics[width=0.39\linewidth]{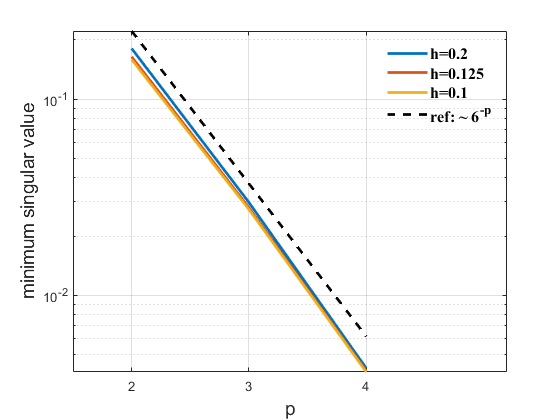}}
      \caption{The maximum/minimum singular value and the spectral condition number 
      of $\mathcal{K}(\mathbf{M}_{1})$ for $d=3$, versus $h$ (at left) and 
      versus $p$ (at right).}
      \label{fig:3d_1_mass}
\end{figure}

Figs.\ref{subfig:1d_1_maxM_h}, \ref{subfig:2d_1_maxM_h} and
      \ref{subfig:3d_1_maxM_h} show that $\sigma_{max}(\mathbf{M}_{1})$ 
      is independent of both $h$ and $p$ for $h>0$ and $p\geq 2$.
In Figs.\ref{fig:1d_1_mass}-\ref{fig:3d_1_mass} we 
      observe that $\sigma_{min}(\mathbf{M}_{1})$ and $\mathcal{K}(\mathbf{M}_{1})$ 
      are also independent of $h$.
For $p-$refinement, the minimum singular value $\sigma_{min}(\mathbf{M}_{1})$ 
      decreases exponentially with $p$.
Specifically, for $d=1$, $\sigma_{min}(\mathbf{M}_{1})$ decreases by about $2^{-p}$;
for $d=2$, it decreases by about $4^{-p}$;
and for $d=3$, it decreases by about $6^{-p}$.
Thus, we obtain the estimations $\sigma_{min}(\mathbf{M}_{1})\thicksim (2d)^{-p}$ 
      and $\mathcal{K}(\mathbf{M}_{1})\thicksim (2d)^{p}$.

\subsection{Sparsity pattern}
It is well known that the number of non-zero entries of the collocation matrix 
      is a measure not only of the memory space 
      required to store the matrix, 
      but also of the computational complexity involved in 
      assembling the collocation matrix and 
      solving the linear system.
In this section, we analyze the sparsity patterns of $\mathbf{A}^{T}\mathbf{A}$ 
      for regularity $k=1$ and $k=p-1$, with three values of degree $p=4,8,12$ and 
      three values of mesh size $h=0.2,0.1,0.067$ in two-dimensional case.
Here, $\mathbf{A}$ denotes the collocation matrix in IGA-L.

\begin{figure}[H]
      \centering 
      \subfigure[h=0.1, p=4, k=1]{
            \includegraphics[width=0.31\linewidth]{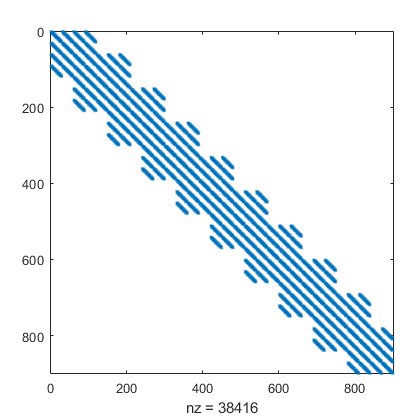}}    
      \subfigure[h=0.1, p=8, k=1]{
            \includegraphics[width=0.31\linewidth]{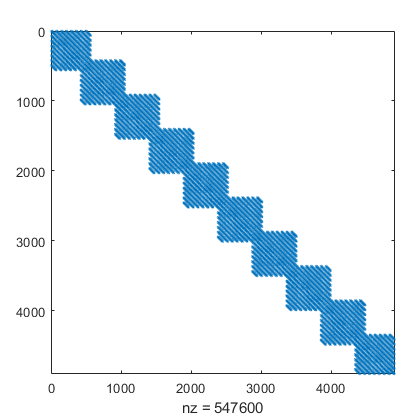}}            
      \subfigure[h=0.1, p=12, k=1]{
            \includegraphics[width=0.31\linewidth]{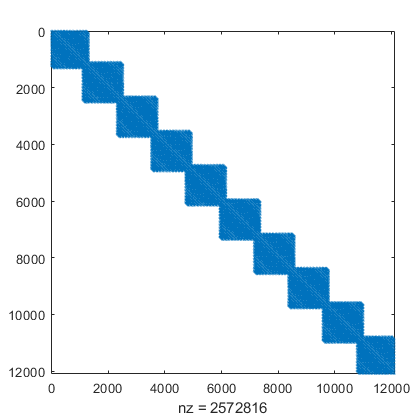}}\\                 
      \subfigure[h=0.1, p=4, k=p-1]{
            \includegraphics[width=0.31\linewidth]{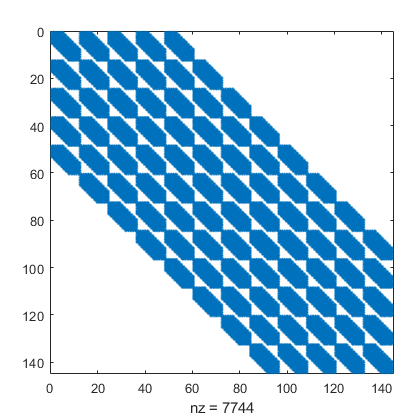}}
      \subfigure[h=0.1, p=8, k=p-1]{
            \includegraphics[width=0.31\linewidth]{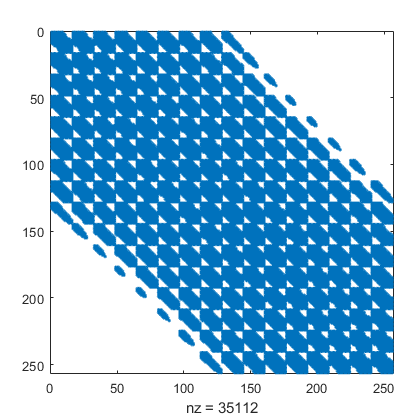}}                
      \subfigure[h=0.1, p=12, k=p-1]{
            \includegraphics[width=0.31\linewidth]{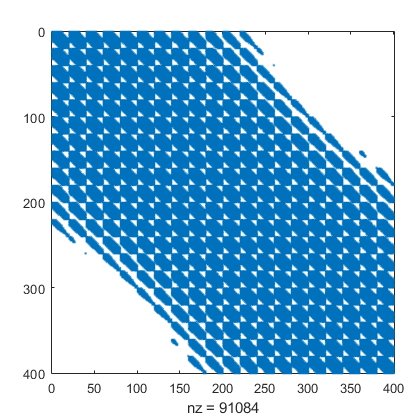}}
      \caption{The sparsity pattern of $\mathbf{A}^{T}\mathbf{A}$ for a fixed mesh size $h=0.1$.}
      \label{fig:sparse_p}
\end{figure}

\begin{figure}[htbp]
      \centering 
      \subfigure[h=0.2, p=8, k=1]{
            \includegraphics[width=0.31\linewidth]{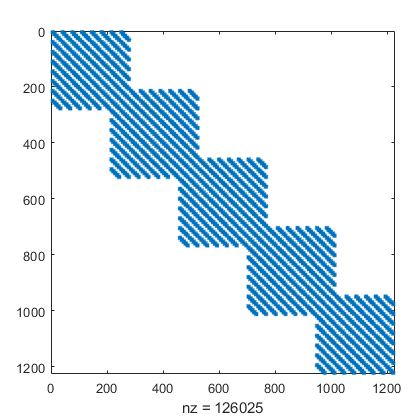}}    
      \subfigure[h=0.1, p=8, k=1]{
            \includegraphics[width=0.31\linewidth]{sparse/h10p8_1.png}}            
      \subfigure[h=0.067, p=8, k=1]{
            \includegraphics[width=0.31\linewidth]{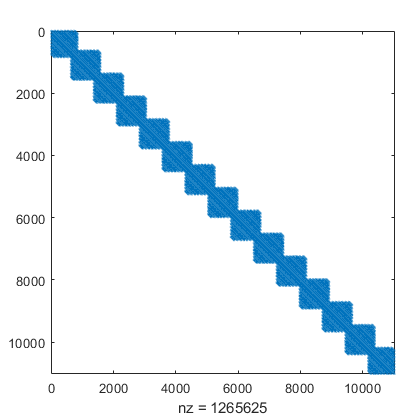}}\\                 
      \subfigure[h=0.2, p=8, k=p-1]{
            \includegraphics[width=0.31\linewidth]{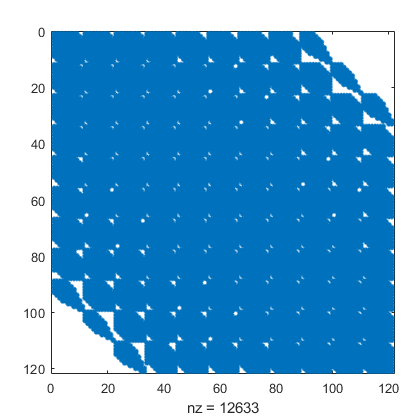}}
      \subfigure[h=0.1, p=8, k=p-1]{
            \includegraphics[width=0.31\linewidth]{sparse/h10p8_p.png}}                
      \subfigure[h=0.067, p=8, k=p-1]{
            \includegraphics[width=0.31\linewidth]{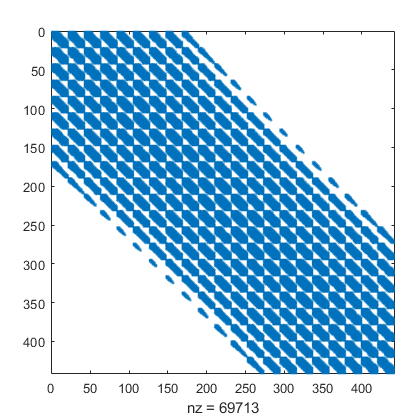}}
      \caption{The sparsity pattern of $\mathbf{A}^{T}\mathbf{A}$ for a fixed degree $p=8$.}
      \label{fig:sparse_h}
\end{figure}

In Fig.\ref{fig:sparse_p}, we illustrate the sparsity patterns of  
      $\mathbf{A}^{T}\mathbf{A}$ for regularity $k=1$ (top), $k=p-1$ (bottom), 
      with three values of degree $p=4$ (left), $p=8$ (center) and 
      $p=12$ (right) for a fixed mesh size $h=0.1$.
We observe that the number of non-zero entries in $\mathbf{A}^{T}\mathbf{A}$ increases 
      as $p$ increases.
Additionally, $\mathbf{A}^{T}\mathbf{A}$ becomes considerably denser as the 
      regularity increases.

Fig.\ref{fig:sparse_h} shows the sparsity patterns of 
      $\mathbf{A}^{T}\mathbf{A}$ for regularity $k=1$ (top), $k=p-1$ (bottom), 
      with three values of mesh size $h=0.2$ (left), $h=0.1$ (center) and 
      $h=0.067$ (right) for a fixed degree $p=8$.
We observe that the number of non-zero entries in $\mathbf{A}^{T}\mathbf{A}$ increases 
      as $h$ decreases.
Compared to IGA-C, the collocation matrix in IGA-L 
      is slightly denser for the same parameters. 
For example, when $h = 0.2$, $p = 8$, and $k=1$, 
      the number of non-zero entries in the 
      collocation matrix of IGA-C is $100,489$, 
      whereas in IGA-L it is $126,025$ \citep{zampieri2024conditioning}.

\subsection{Collocation point}
In this section, we analyze the effect of the number and the distribution 
      of the collocation points on the spectral properties of the collocation 
      matrix.
The degree $p$ is fixed at $8$, and the mesh size is set to $h=0.05$.
We conduct experiments in one-dimensional case, 
      with the number of collocation points 
      initially set to be equal to the $dof$.

\begin{figure}[htbp]
            \centering 
            \subfigure[Condition number vs. m (k=p-1)]{
                  \label{subfig:numcol_cond_p}
                  \includegraphics[width=0.46\linewidth]{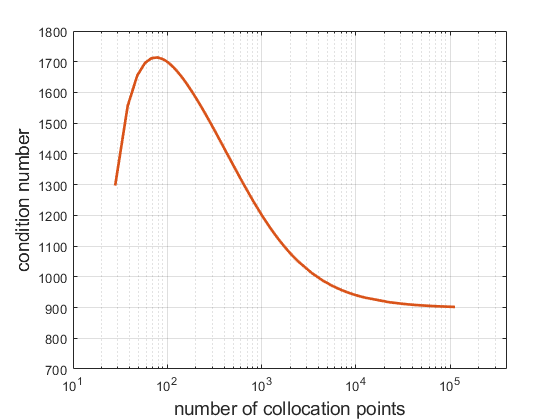}}    
            \subfigure[Condition number vs. m (k=1)]{
                  \label{subfig:numcol_cond_1}
                  \includegraphics[width=0.46\linewidth]{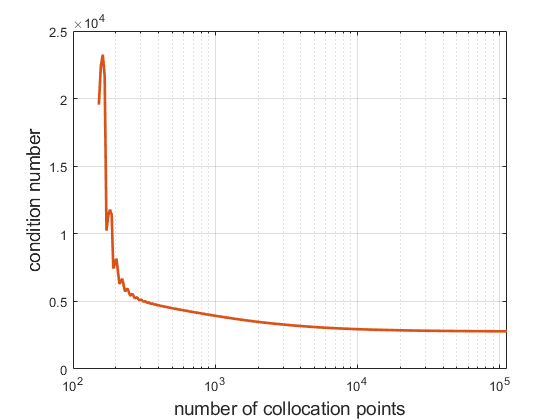}}\\             
            \subfigure[Maximum singular value vs. m (k=p-1)]{
                  \label{subfig:numcol_max_p}
                  \includegraphics[width=0.46\linewidth]{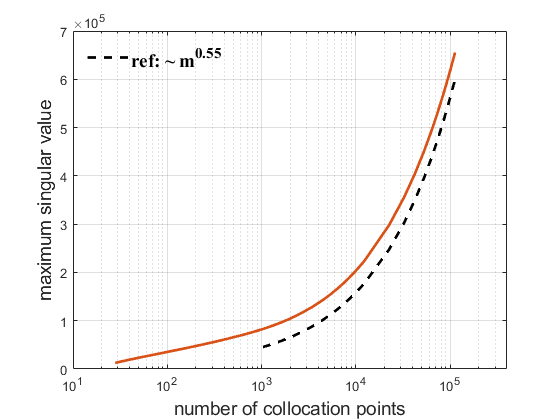}}                
            \subfigure[Maximum singular value vs. m (k=1)]{
                  \label{subfig:numcol_max_1}
                  \includegraphics[width=0.46\linewidth]{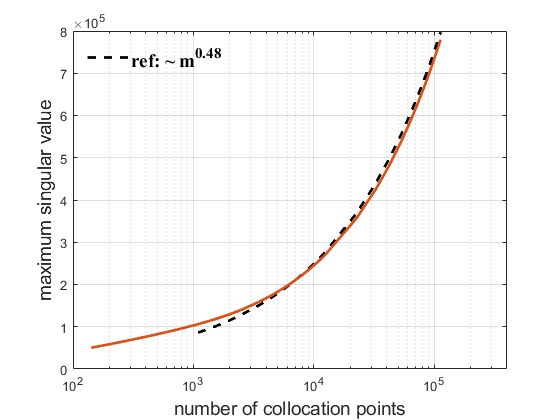}}\\
            \subfigure[Minimum singular value vs. m (k=p-1)]{
                  \label{subfig:numcol_min_p}
                  \includegraphics[width=0.46\linewidth]{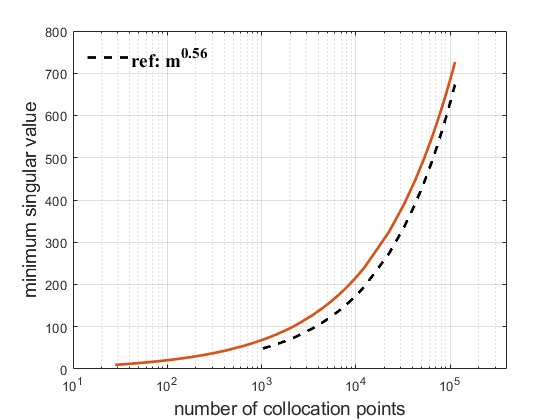}}                
            \subfigure[Minimum singular value vs. m (k=1)]{
                  \label{subfig:numcol_min_1}
                  \includegraphics[width=0.46\linewidth]{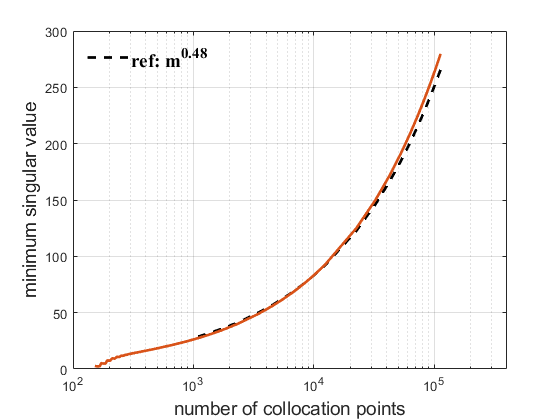}}
            \caption{The maximum/minimum singular value and the 
            spectral condition number of the collocation matrix for different 
            number of collocation points.}
            \label{fig:num_collopoints}
      \end{figure}

Fig.\ref{fig:num_collopoints} shows the condition number (top), 
      maximum singular value (center) and minimum singular value (bottom) 
      as the number of collocation points $m$ increases 
      for $k=p-1$ (left) and $k=1$ (right).
In Figs.\ref{subfig:numcol_cond_p} and \ref{subfig:numcol_cond_1} we 
      observe that as the number of collocation points increases, 
      $\mathcal{K}(A)$ initially increases, then decreases, and 
      finally stabilizes.
When the number of collocation points significantly exceeds the number of $dof$, 
      $\mathcal{K}(A)$ can be much lower than that in IGA-C.
For $k=1$ regularity, the condition number decreases more rapidly than for $k=p-1$ 
      as the number of collocation points increases, 
      suggesting that in the IGA-L method with $k=1$ regularity, 
      using slightly more collocation points than the $dof$ 
      can lead to a better-conditioned linear system compared to IGA-C.

\begin{figure}[htbp]
      \centering 
      \subfigure[Condition number vs. h (SC)]{
            \label{subfig:sc_cond_h}
            \includegraphics[width=0.47\linewidth]{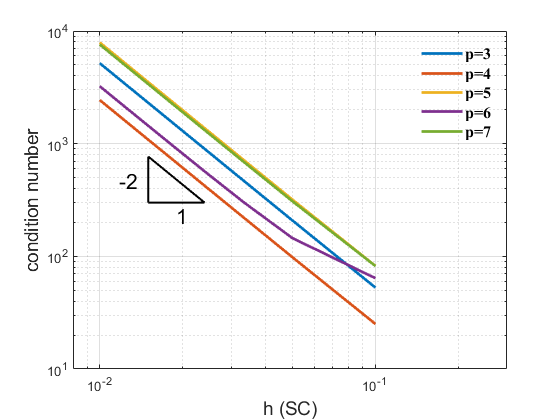}}    
      \subfigure[Condition number vs. h (CG)]{
            \label{subfig:cg_cond_h}
            \includegraphics[width=0.47\linewidth]{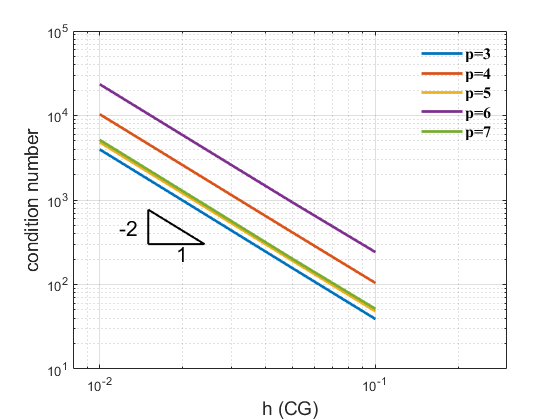}} \\                         
      \subfigure[Condition number vs. $p$ (SC)]{
            \label{subfig:sc_cond_p}
            \includegraphics[width=0.47\linewidth]{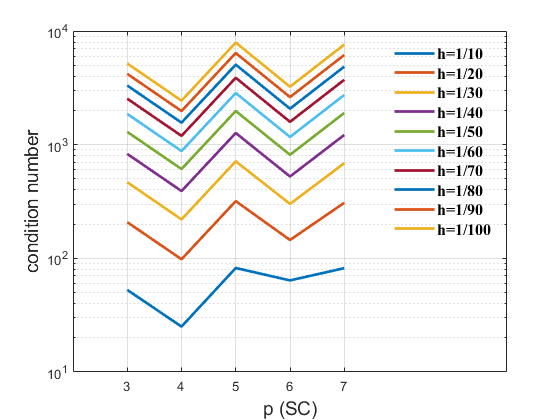}}
      \subfigure[Condition number vs. $p$ (CG)]{
            \label{subfig:cg_cond_p}
            \includegraphics[width=0.47\linewidth]{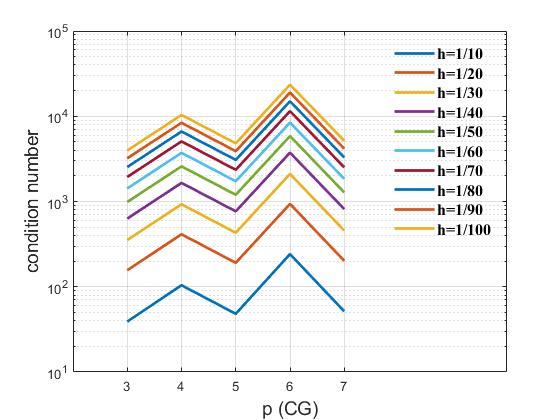}} \\
      \subfigure[Maximum singular value vs. $h$ (SC)]{
            \label{subfig:sc_max_p}
            \includegraphics[width=0.47\linewidth]{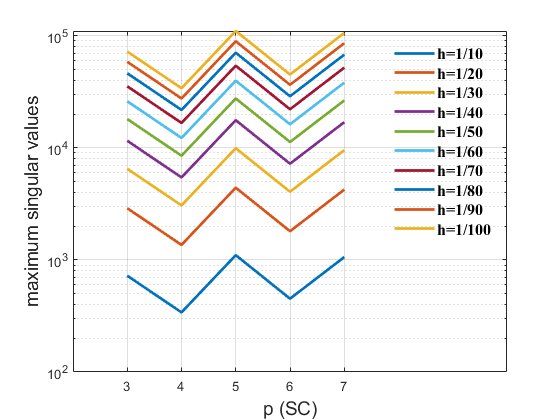}}                       
      \subfigure[Maximum singular value vs. $p$ (CG)]{
            \label{subfig:cg_max_p}
            \includegraphics[width=0.47\linewidth]{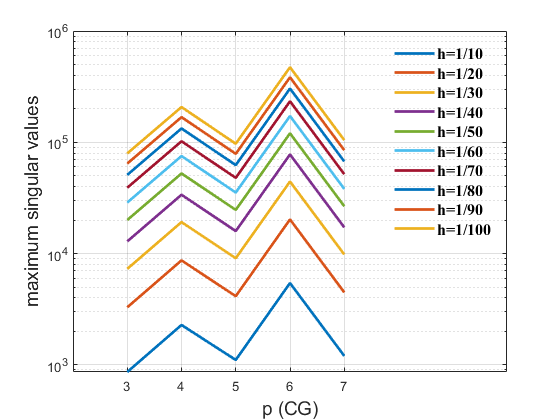}}
      \caption{The spectral condition number and maximum singular value of 
       the collocation matrix formed by SC points and CG points.}
      \label{fig:SC_CGpoints}
\end{figure}

\begin{figure}[htbp]
      \centering        
      \subfigure[SC points ($p=4$)]{
            \includegraphics[width=0.92\linewidth]{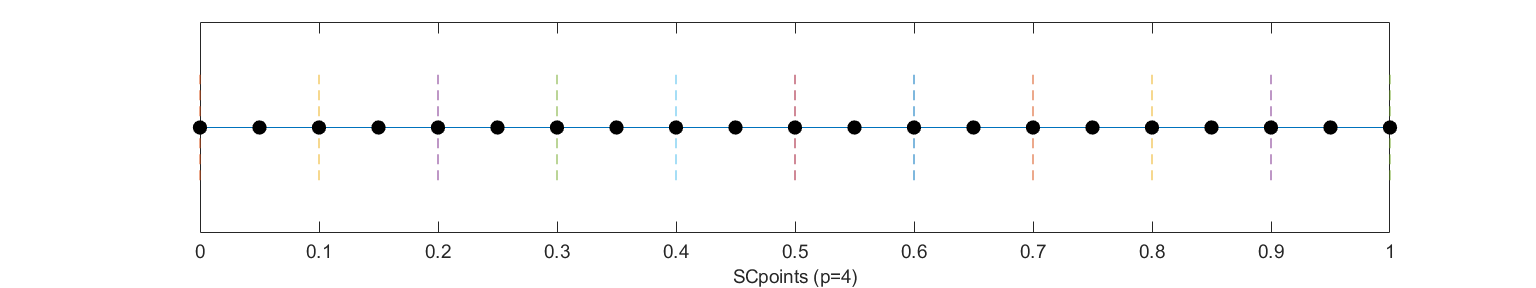}}\\                
      \subfigure[CG points ($p=4$)]{
            \includegraphics[width=0.92\linewidth]{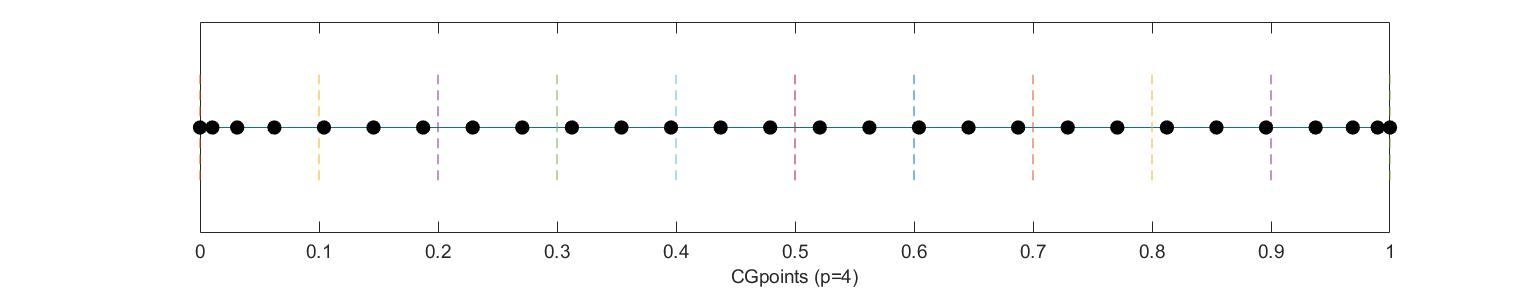}}
      \caption{Distribution of SC points (top) and CG points (bottom) when $p=4$.}
      \label{fig:distribution_SC_CGpoints}
\end{figure}

From Figs.\ref{subfig:numcol_max_p} and \ref{subfig:numcol_max_1} we 
      observe that as the number of collocation points increases, 
      $\sigma_{max}(\mathbf{A}_{p-1})\thicksim m^{0.55}$ and 
      $\sigma_{max}(\mathbf{A}_{1})\thicksim m^{0.48}$, 
      indicating that the maximum singular value for $k=p-1$ regularity 
      increases more rapidly than for $k=1$.
Similarly, Figs.\ref{subfig:numcol_max_p} and \ref{subfig:numcol_max_1} show that 
      $\sigma_{min}(\mathbf{A}_{p-1})\thicksim m^{0.55}$ and 
      $\sigma_{min}(\mathbf{A}_{1})\thicksim m^{0.48}$.
Fig.\ref{fig:SC_CGpoints} shows the condition number of 
      the collocation matrix formed by SC points and CG points 
      versus $h$ and $p$.
From Figs.\ref{subfig:sc_cond_h} and \ref{subfig:cg_cond_h} we observe that 
      the relationship between the condition number and the mesh size $h$ 
      is similar to the previous results obtained using Greville points.
Figs.\ref{subfig:sc_cond_p} and \ref{subfig:cg_cond_p} show the 
      fluctuation of the condition numbers related to both SC points and 
      CG points for $p-$refinement.
The condition numbers of the collocation matrix formed by SC points 
     in odd-degree cases are lower than those in even-degree cases.
On the contrary, for CG points, the spectral properties in even-degree cases 
      are superior to those in odd-degree cases.
As shown in Figs.\ref{subfig:sc_max_p} and \ref{subfig:cg_max_p}, 
      we find that the variations in condition numbers 
      corresponding to different distributions of collocation points 
      are primarily influenced by the maximum singular values.
From Figs.\ref{fig:dis_collocation} and \ref{fig:distribution_SC_CGpoints}, 
      we observe that the SC points in odd-degree cases 
      are more densely distributed near the boundary 
      than those in even-degree cases. 
Similarly, the CG points show a denser distribution 
      near the boundary for even-degrees than for odd-degrees.
The dense distribution of the collocation points 
      increases the $L_{2}$ norm of the column vectors 
      corresponding to the basis functions near the boundary,  
      resulting in an uneven distribution of the column vectors.  
Consequently, the maximum singular value and the condition number increase.
A similar explanation can account for the fluctuation of the condition number 
      for CG points with different degrees.

\section{Conclusion}
\label{sec:con}
In this paper, we conduct a systematic analysis 
      of the spectral properties of the collocation matrix and mass matrix 
      in IGA-L arising from Poisson problem with 
      homogeneous Dirichlet boundary conditions.
We focus on the relationship between the condition number, maximum and minimum singular 
      value and the IGA-L discretization parameters $h$, $p$, $k$  
      and the number and distribution of the collocation points.
Through numerical investigation, we present theoretical estimations of the condition number, 
      maximum, and minimum singular value of the collocation matrix and mass matrix,  
      filling a gap in the theoretical research in this area.
The condition number of the IGA-L collocation matrix for $k=p-1$ regularity  
      exhibits different behavior in different regions of the $(p,h)$ plane.
For small $h$, the condition number behaves as $h^{-2}p^{2}$,  
      while for large $h$, it becomes independent of $h$ 
      and grows exponentially with $p$.
However, the growth rate is slower than that of the stiffness matrix 
      in IGA-G under the same conditions.
Therefore, the linear system in IGA-L is better conditioned compared to 
      IGA-G, especially as $p$ increases.
For the mass matrix, its condition number is independent of $h$ 
      and grows exponentially with $p$, 
      with a faster growth rate for $k=1$ than for $k=p-1$.

Additionally, we observe that for IGA-L with $k=1$ regularity, 
      using slightly more collocation points than the number of $dof$ 
      can lead to a better-conditioned linear system.
Furthermore, the collocation matrix for $k=1$ contains significantly 
      fewer non-zero entries than for $k=p-1$,
      reducing the computational cost of assembling the matrix 
      and solving the linear system.

Future work will focus on enhancing the spectral properties of the collocation matrix 
      by optimizing the distribution and number of collocation points, 
      as well as determining appropriate ranges for mesh size $h$ and degree $p$.
Additionally, we will also consider the extension of this work to 
      advection-diffusion-reaction problems.  

\section*{References}

\bibliography{mybibfile}

\begin{thebibliography}{10}
\expandafter\ifx\csname url\endcsname\relax
  \def\url#1{\texttt{#1}}\fi
\expandafter\ifx\csname urlprefix\endcsname\relax\def\urlprefix{URL }\fi
\expandafter\ifx\csname href\endcsname\relax
  \def\href#1#2{#2} \def\path#1{#1}\fi

\bibitem{Hughes2005IGA}
T.~J.~R. Hughes, J.~A. Cottrell, Y.~Bazilevs, Isogeometric analysis: Cad,
  finite elements, nurbs, exact geometry and mesh refinement, Computer Methods
  in Applied Mechanics and Engineering 194~(39-41) (2005) 4135--4195.
\newblock \href {http://dx.doi.org/10.1016/j.cma.2004.10.008}
  {\path{doi:10.1016/j.cma.2004.10.008}}.

\bibitem{bazilevs2006isogeometric}
Y.~Bazilevs, L.~Beirao~da Veiga, J.~A. Cottrell, T.~J. Hughes, G.~Sangalli,
  Isogeometric analysis: approximation, stability and error estimates for
  h-refined meshes, Mathematical Models and Methods in Applied Sciences 16~(07)
  (2006) 1031--1090.

\bibitem{cottrell2009isogeometric}
J.~A. Cottrell, T.~J. Hughes, Y.~Bazilevs, Isogeometric analysis: toward
  integration of CAD and FEA, John Wiley \& Sons, 2009.

\bibitem{da2014mathematical}
L.~B. Da~Veiga, A.~Buffa, G.~Sangalli, R.~V{\'a}zquez, Mathematical analysis of
  variational isogeometric methods, Acta Numerica 23 (2014) 157--287.

\bibitem{hughes2008duality}
T.~J. Hughes, A.~Reali, G.~Sangalli, Duality and unified analysis of discrete
  approximations in structural dynamics and wave propagation: comparison of
  p-method finite elements with k-method nurbs, Computer methods in applied
  mechanics and engineering 197~(49-50) (2008) 4104--4124.

\bibitem{Auricchio2011IGAC}
F.~Auricchio, L.~B. Da~Veiga, T.~J.~R. Hughes, A.~Reali, G.~Sangalli,
  Isogeometric collocation methods, Mathematical Models and Methods in Applied
  Sciences 20~(11) (2011) 2075--2107.
\newblock \href {http://dx.doi.org/10.1142/s0218202510004878}
  {\path{doi:10.1142/s0218202510004878}}.

\bibitem{auricchio2012isogeometric}
F.~Auricchio, L.~B. Da~Veiga, T.~J. Hughes, A.~Reali, G.~Sangalli, Isogeometric
  collocation for elastostatics and explicit dynamics, Computer methods in
  applied mechanics and engineering 249 (2012) 2--14.

\bibitem{manni2015isogeometric}
C.~Manni, A.~Reali, H.~Speleers, Isogeometric collocation methods with
  generalized b-splines, Computers \& Mathematics with Applications 70~(7)
  (2015) 1659--1675.

\bibitem{anitescu}
C.~Anitescu, Y.~Jia, Y.~J. Zhang, T.~Rabczuk, An isogeometric collocation
  method using superconvergent points, Computer Methods in Applied Mechanics
  and Engineering 284 (2015) 1073--1097.
\newblock \href {http://dx.doi.org/10.1016/j.cma.2014.11.038}
  {\path{doi:10.1016/j.cma.2014.11.038}}.

\bibitem{LinIGA-L}
H.~Lin, Y.~Xiong, X.~Wang, Q.~Hu, J.~Ren, Isogeometric least-squares
  collocation method with consistency and convergence analysis, JOURNAL OF
  SYSTEMS SCIENCE \& COMPLEXITY 33~(5) (2020) 1656--1693.
\newblock \href {http://dx.doi.org/10.1007/s11424-020-9052-9}
  {\path{doi:10.1007/s11424-020-9052-9}}.

\bibitem{demko}
S.~DEMKO, On the existence of interpolating projections onto spline spaces,
  Journal OF Approximation Theory 43~(2) (1985) 151--156.
\newblock \href {http://dx.doi.org/10.1016/0021-9045(85)90123-6}
  {\path{doi:10.1016/0021-9045(85)90123-6}}.

\bibitem{Gomez}
H.~Gomez, L.~De~Lorenzis, The variational collocation method, Computer Methods
  in Applied Mechanics and Engineering 309 (2016) 152--181.
\newblock \href {http://dx.doi.org/10.1016/j.cma.2016.06.003}
  {\path{doi:10.1016/j.cma.2016.06.003}}.

\bibitem{manni2011generalized}
C.~Manni, F.~Pelosi, M.~L. Sampoli, Generalized b-splines as a tool in
  isogeometric analysis, Computer Methods in Applied Mechanics and Engineering
  200~(5-8) (2011) 867--881.

\bibitem{bazilevs2010isogeometric}
Y.~Bazilevs, V.~M. Calo, J.~A. Cottrell, J.~A. Evans, T.~J.~R. Hughes,
  S.~Lipton, M.~A. Scott, T.~W. Sederberg, Isogeometric analysis using
  t-splines, Computer methods in applied mechanics and engineering 199~(5-8)
  (2010) 229--263.

\bibitem{bathe2006finite}
K.-J. Bathe, Finite element procedures, Klaus-Jurgen Bathe, 2006.

\bibitem{gahalaut2012condition}
K.~Gahalaut, S.~Tomar, Condition number estimates for matrices arising in the
  isogeometric discretizations, Johann Radon Institut (RICAM), 2012.

\bibitem{garoni2014spectrum}
C.~Garoni, C.~Manni, F.~Pelosi, S.~Serra-Capizzano, H.~Speleers, On the
  spectrum of stiffness matrices arising from isogeometric analysis, Numerische
  Mathematik 127 (2014) 751--799.

\bibitem{eisentrager2020condition}
S.~Eisentr{\"a}ger, E.~Atroshchenko, R.~Makvandi, On the condition number of
  high order finite element methods: Influence of p-refinement and mesh
  distortion, Computers \& Mathematics with Applications 80~(11) (2020)
  2289--2339.

\bibitem{szabo2021finite}
B.~Szab{\'o}, I.~Babu{\v{s}}ka, Finite element analysis: Method, verification
  and validation.

\bibitem{canuto2007spectral}
C.~Canuto, M.~Y. Hussaini, A.~Quarteroni, T.~A. Zang, Spectral methods:
  evolution to complex geometries and applications to fluid dynamics, Springer
  Science \& Business Media, 2007.

\bibitem{solin2003higher}
P.~Solin, K.~Segeth, I.~Dolezel, Higher-order finite element methods, Chapman
  and Hall/CRC, 2003.

\bibitem{demkowicz2006computing}
L.~Demkowicz, Computing with hp-adaptive finite elements: volume 1 one and two
  dimensional elliptic and Maxwell problems, Chapman and Hall/CRC, 2006.

\bibitem{melenk2002condition}
J.~M. Melenk, On condition numbers in hp-fem with gauss--lobatto-based shape
  functions, Journal of Computational and Applied Mathematics 139~(1) (2002)
  21--48.

\bibitem{gervasio2020computational}
P.~Gervasio, L.~Ded{\`e}, O.~Chanon, A.~Quarteroni, A computational comparison
  between isogeometric analysis and spectral element methods: accuracy and
  spectral properties, Journal of Scientific Computing 83 (2020) 1--45.

\bibitem{zampieri2024conditioning}
E.~Zampieri, L.~F. Pavarino, Conditioning and spectral properties of
  isogeometric collocation matrices for acoustic wave problems, Advances in
  Computational Mathematics 50~(2) (2024) 16.

\bibitem{piegl2012nurbs}
L.~Piegl, W.~Tiller, The NURBS book, Springer Science \& Business Media, 2012.

\end{thebibliography}

\end{document}